\documentclass[12pt, singlespace]{amsart}

\usepackage{amsmath,amssymb,amsthm,epsfig }
\usepackage{enumerate}
\usepackage[dvips]{color}
\usepackage{bbm}
\usepackage[colorinlistoftodos]{todonotes}
\usepackage{version}
\usepackage[margin=20truemm]{geometry}
\usepackage{mathtools}
\usepackage[all, cmtip]{xy}
\SelectTips{cm}{10}

\newtagform{alph}[]()
\newtagform{Alph}[]()

\DeclareFontFamily{OT1}{rsfs}{}
\DeclareFontShape{OT1}{rsfs}{n}{it}{<-> rsfs10}{}
\DeclareMathAlphabet{\curly}{OT1}{rsfs}{n}{it}

\newtheorem{Thm}{Theorem}[section]
\newtheorem{lem}[Thm]{Lemma}

\newtheorem{cor}[Thm]{Corollary}

\newtheorem{prop}[Thm]{Proposition}

\newtheorem{``Conj"}[Thm]{``Conjecture"}

\newtheorem{Claim}{Claim}
\newtheorem{Property}{Property}

\theoremstyle{remark}
\newtheorem{Rem}[Thm]{Remark}
\newtheorem{ex}[Thm]{Example}
\theoremstyle{definition}
\newtheorem{defn}[Thm]{Definition}

\newtheorem{Step}{Step}

\newtheorem{ntt}{Notation}

\newcommand{\Spec}{\mathop{\mathrm{Spec}}\nolimits}

\newcommand{\gr}{\mathop{\mathrm{gr}}\nolimits}

\newcommand{\bfit}[1]{\textbf{\textit{#1}}}

\def\k{\mathbbm{k}}
\def\R{\mathbb{R}}

\def\C{\mathbb{C}}
\def\G{\mathbb{G}}
\def\Z{\mathbb{Z}}
\def\O{\mathcal{O}}

\def\Q{\mathbb{Q}}

\def\A{\mathbb{A}}

\def\X{\mathcal{X}}
\def\Y{\mathcal{Y}}

\def\F{\mathcal{F}}

\excludeversion{memo}

\begin{document}

\title[Canonical degenerations over toroidal base]
{Stability theory over toroidal 
or Novikov-type 
base and canonical modifications}
\author{Yuji Odaka}

\maketitle
\thispagestyle{empty}

\begin{abstract}
We set up a generalization of ubiquitous 
one-parameter families in algebraic geometry and their 
use for stability theories 
(\cite{GIT, HL, AHLH}) 
to families over toric varieties and their analytic analogues. 
The language allows us to 
reformulate degenerations of ``irrational" direction in the literature 
as {\it canonical} objects in a unified manner. 

Accordingly, 
we generalize the 
(semi)stable reduction-type theorem 
for $\Theta$-stratification 
in \cite{AHLH} of Langton type 
to our higher rank setup. We also establish a complex analytic 
analogue of the results. As an infinitesimal analogue of toric spectrum, 
we also use Novikov type rings, as they provide greater 
canonicity, 
but their use can be 
avoided logically for readers who prefer not to use such rings. 
Our other papers discuss applications of 
our framework and the main theorems to 
problems around K-stability and K\"ahler 
geometry. 
\end{abstract}

\section{Introduction}\label{sec:intro}


In the classical geometric invariant theory \cite{GIT}, 
the so-called 
Hilbert-Mumford criterion for stability of an algebro-geometric 
object, say ${\bf X}$, with respect to an action of a reductive algebraic group 
$G$, shows that it is enough to analyze degenerations of ${\bf X}$ 
along actions of $\G_m$ which is inside $G$ ($1$-parameter subgroups). 
The resulting families are isotrivially degenerating 
families of ${\bf X}$ over $\A^1$. 
In the case when ${\bf X}$ is a polarized 
projective variety, 
they are characterized as {\it test configurations} by 
Donaldson \cite{Don02} (cf., also \cite{Don05}). 

More precisely, for the Hilbert-Mumford criterion, 
one associates a GIT weight to 
each $1$-parameter subgroup and their positivity is equivalent to 
the stability. 
Alternatively, it is well known that one can consider maximal tori $T$ of $G$, instead of 
$1$-parameter subgroups, 
and encode all GIT weights data 
within $T$ 
in the {\it weight polytope}. 
The corresponding natural extension of test configurations or 
the above type isotrivial degenerations, replacing $\G_m$ by 
a maximal torus, is isotrivial 
family  
over $T$ and its partial equivariant compactification i.e., 
toric varieties. 
A simple toy model for the principal objects 
in this paper is an isotrivial degeneration over a higher-dimensional toric variety. 

More generally, we introduce 
more extensive 
algebraic (partly stack-theoretic) 
generalities with such perspective, which leads to  introduce e.g., 
isotrivial degenerating families over toroidal bases, 
Novikov type valuation rings 
of high rank, and their analytic analogues. 
In particular, focusing on 
the context of stability and moduli theories of 
varieties, we introduce the  notion of {\it generalized test configurations} and 
their limit versions. Then, together with a companion paper, 
we show that the generalized notion is useful 
to treat the following types of degenerations in ``irrational direction" 
in the literature, in a more canonical and unified manner: 
\begin{itemize}
    \item klt singularities to irregular Ricci-flat cones (\cite{DSII, Li}), 
    \item K\"ahler-Ricci flow on Fano manifolds converging to K\"ahler-Ricci solitons (\cite{CSW, HanLi}), 
    \item negative weight deformations to asymptotically conical Calabi-Yau manifolds (\cite{CH}). 
\end{itemize}

In our other papers (cf., \cite{Od24a, Od24c}), 
we discuss and explore these examples one by one, but 
this paper we focus on the general framework. 
In particular, our \S \ref{sec:3} 
gives some generalization of 
the $\Theta$-(semi)stability 
theory and semistable reduction theorem (\cite{HL} and \cite{AHLH} 
respectively), as a core theorem in this paper. 
Here, recall $\Theta:=[\A^1/\G_m]$ over a field $k$ 
is the quotient stack whose coarse moduli is just one point which is 
used as a source of map to express equivariant isotrivial families over 
$\A^1$. \cite{HL, AHLH} introduces the notion of 
$\Theta$-strata and $\Theta$-stratification on 
(highly non-separated) algebraic stacks, as stacky analogues of 
the Hesselink-Kempf-Kirwan-Ness stratification  \cite{Kirwan, Ness} and 
the Harder-Narasimhan filtration of coherent sheaves \cite{HN}. 
We call a natural 
higher rank($=:r$) version of it as higher $\Theta$-strata 
by roughly considering degenerations along $\Theta^r=[\A^r/T]$ 
with $T:=\G_m^r$ (temporarily). 
One can also replace $\A^r$ by other affine toric variety $U_\tau$ for a regular rational polyhedral cone 
$\tau$. Note that due to the regularity assumption, which we assume 
throughout the paper, 
$U_\tau\simeq \A^r$. 
\footnote{We believe this assumption is redundant but to use earlier 
papers which assume this, we put it for simplicity. In practice, 
we are normally interested in sufficiently small $\tau$, so it does not 
change the practicality.} 

To give a feel of our main statements on such 
higher $\Theta$-stratum, 
here we show its quick brief outlook, 
although we refer to later 
Theorem \ref{gAHLH}, \ref{gAHLH.R}, Corollary \ref{gAHLHc} 
for the precise statements 
and generalizations. Roughly put, 
it shows the existence of ``improvement/modification" 
of one parameter family, which 
avoids some ``too unstable" locus in a fairly general setup. 

For below, we consider a 
quotient stack  $\mathcal{M}=[M/G]$ over $\k$, 
where $M$ is an algebraic $\k$-scheme and 
$G$ is a linear algebraic $\k$-group (e.g., ${\rm GL}_n$ 
or algebraic torus). 
In that case, higher $\Theta$-stratum $\mathcal{Z}^+=[Z/G]$ 
with the closed subscheme $Z\subset M$ and a rational regular 
polyhedral cone $\tau\subset N\otimes \R$ 
roughly means for each point $p\in Z(\k)$, 
there is its isotrivial degeneration (the generalized 
Harder-Narasimhan filtration) along the affine 
toric variety $U_\tau\supset T$ with a homomorphism 
$T\to G$ and $Z^+$ is the maximal locus 
where this structure varies uniformly. 
Note that our $r$ is not necessarily the rank of $G$. 

\begin{Thm}[Higher $\Theta$-semistable reduction: cf., Corollary \ref{gAHLHc}]\label{intro.Mthm}
Assume $\k$ is algebraically closed field of characteristic $0$ (e.g., 
$\k=\C$). 
Consider a 
quotient 
stack $\mathcal{M}$ over 
a field $\k$ 
and its 
higher $\Theta$-strata $\mathcal{Z}^+\subset \mathcal{M}$ for the cone 
(see Definition \ref{higher.theta.strata}). 

Take any $\k$-morphism $f\colon \Delta=\Spec \k[[t]]\to \mathcal{M}$ 
such that its closed point 
maps into $\mathcal{Z}^+\subset  \mathcal{M}$ 
while the generic point 
maps outside $\mathcal{Z}^+$. 
We denote the restriction of $f$ to $\Spec(K)$ 
as $f^o$. 

Then, for $\Spec\k[[t^{1/d}]]$ for some $d\in \Z_{>0}$, 
there exists a 
modification of $f$ as follows: 

$$f_{{\rm new}} \colon \Spec\k[[t^{1/d}]]
\to \mathcal{M},$$
which extends the base change of $f^o$ to $\k((t^{\frac{1}{d}}))$ but 
sends the closed point $c_{R'}$ to {\bf outside} $\mathcal{Z}^+$. 
\end{Thm}
If there are finitely many strata $\{\mathcal{Z}_i^+\}_i$ inside  $\mathcal{M}$ 
which form what we call  
{\it higher $\Theta$-stratification} (Definition \ref{Theta.strat}), extending the case $\#\{i\}=1$ above, 
we can iteratively apply the above modification 
(see Theorem \ref{theta.stratification.
stablereduction}) 
to obtain the semistable reduction for the complement of $\mathcal{Z}_i^+$s, 
which often admit moduli interpretations. 
Indeed, in the moduli contexts, 
$\mathcal{Z}_i^+$s often appear as unstable locus 
in the spirit of \cite{GIT} so that their complement 
is nothing but the locus which parametrizes semistable objects. 

\vspace{2mm}
Here are some technical comments in order.
For Theorem \ref{intro.Mthm} mentioned above,
we also refine the statements in the main body of this text
from at least the following four different aspects.

The first aspect is about the new base $R'$. 
There are many choices of the pair of DVR 
$R'$ and $f_{\rm new}$ {\it as a family} 
one can take in our construction, 
though $f_{\rm new}(c_{R'})$ is more canonical. 
Indeed, they roughly correspond to choices of 
integral points in a certain rational polyhedral cone so that there is no canonical choice and there are infinite possibilities. 
Nevertheless, 
if we fix a possibly irrational vector $\xi(\neq 0)\in \tau$ and 
allow $R'$ to be non-discrete valuation ring in general, 
there is more canonical choice of $R'$ as a bigger valuation  ring for each $\xi$ and we expect its importance for 
future study. That is the Novikov type ring we mentioned above, 
as we introduce and 
discuss in \S 2. 
Secondly, we also prove its complex analytic version (Corollary \ref{gAHLH.R}). 

We call this replacement 
of $f$ {\it higher $\Theta$-semistable 
reduction} and call its limit 
{\it higher $\Theta$-semistable limit}. 
Our result partially recovers 
\cite[\S 6, Theorem 6.3, Theorem 6.5]{AHLH} by a different approach. 

Thirdly, in the case when $\mathcal{M}$ is a quotient stack by an 
algebraic torus, we prove {\it certain uniqueness} of the above 
modification. 
See \eqref{uniqueext} of Theorem \ref{gAHLH} for the details. 
\footnote{This is beneficial for the study of bubbling theory \cite{Od24c}.}

Lastly, but almost obviously, if there are certain stratifications of $\mathcal{M}$ 
consists of higher $\Theta$-strata, we can simply 
repeat the above modification 
to give a stronger result (Theorem \ref{theta.stratification.
stablereduction}) as it was in \cite{AHLH}. 

\vspace{2mm}

The framework and theorems in this paper
are applied in our other papers (cf., 
\cite{Od24a, Od24c}), which are mostly in the 
context of K\"ahler geometry and algebraic 
theory related to K-stability. For instance, 
\cite{Od24a} uses the above Theorem 
\ref{intro.Mthm} 
to establish properness of moduli of 
Calabi-Yau cones and \cite{Od24c} 
uses it to provide an 
algebro-geometric version of 
construction of certain 
{\it bubbling} K\"ahler metrics (its arXiv v2 
\S 2.3 also includes 
application to properness of 
moduli of K\"ahler-Ricci solitons). 

\subsubsection*{Notes added} 
About $8$ and a half months after the first version of this paper, 
we saw some developments i.e., \cite{BHLINK} appeared, but with 
a different motivation. 
At the end of the paper, put in our version 3 on arXiv, we discuss the  overlap, interesting relations and certain 
differences between this paper and \cite{BHLINK}. 

\section{Generalities on multi-parameter families}\label{sec:gentc}

This section is of preparatory nature and shows our initial motivation, 
which is closely related to various asymptotics of K\"ahler metrics 
and some variants of K-stability theory. 
For the readers who are only interested in the general stacky 
statements, we can suggest to look only Lemma \ref{lem:MGamma} to 
Proposition 2.3, and then skip to the next section. 

In this section, we generalize 
one parameter families and degenerations, which appear 
in the theory of K-stability or its local analogues for klt singularities (cf., \cite{Li} and subsequent 
developments), and the stacky stability theories in \cite{HL, AHLH}, 
to those over toroidal base of often higher dimension. 
Some germs of our theory can be seen in \cite{Spi, Teissier, DSII, LX, HanLi}, for instance. 
We observe that our language allows re-explanation in a unified manner of those 
existing K-stability theories. 

Firstly, we introduce a notion of {\it generalized test configurations} as certain type of 
multi-parameter isotrivial family over toroidal bases and their analogues, 
which unifies various families in the field of K-stability and K\"ahler geometry, such as 
\begin{enumerate}
\item 
test configurations and their affine versions (\cite{Don02, CS, CS2}) cf., Example \ref{ex.tc}, 
\item \label{ex:nwd}
negative weight deformation in the sense of \cite[Definition 1.7]{CH} cf., Example \ref{ex.nwd} and 
\item \label{ex:2step}
each step of the recent $2$-step degeneration theory (originated in \cite{DSII}) 
\begin{itemize}
\item 
e.g., 
the canonical degeneration of klt singularities to its {\it possibly irregular}
K-semistable Fano cone via the 
volume-minimizing valuation of klt singularities (\cite{Li}) 
in Example \ref{loc.ex}, 
\item $\R$-test configurations in \cite{CSW, HanLi} cf., Example \ref{ex.r-tc},  
\item (the $2$nd steps to the polystable replacements (\cite{DSII, LWX}) as Examples 
\ref{loc.ex2}), 
\end{itemize}
\item (or more generally, degeneration of GIT semistable objects to GIT polystable objects 
(\cite{GIT}) cf., Remark \ref{ex.GIT}). 
\end{enumerate}
Further, as we will later explain, the notion allows us to treat the above examples of degenerations 
\eqref{ex:nwd} and \eqref{ex:2step} in a (more) canonical manner, to remove the broad 
ambiguity of choices of $1$-parameter 
subgroups, except for those 
with parentheses. 

\subsection{Setup of the bases}

Since our bases will be often toric geometric, first we set the usual 
notation for toric geometry. 

\begin{ntt}\label{not.hr}
We fix a base ring $A$ and a free abelian group (lattice) 
$N$ of rank $r$ and its strictly convex full-dimensional 
rational polyhedral cone $\tau\subset N_{\R}:=N\otimes_{\Z}\R$. 
The dual lattice of $N$ is denoted as $M:={\rm Hom}(N,\Z)$ 
and we set $M_{\R}:=M\otimes \R$. 
We define the algebraic split torus $T=\G_m\otimes_{\Z} N$ as 
${\rm Spec} A[M]$, where $A[M]$ means the group ring over $A$. 
If we regard $M$ as the character lattice of $T$, 
we denote by $\chi^{m}$ the character of $T$ 
which corresponds to each $m\in M$. 
The affine toric partial compactification of $T$ for $\tau$ 
is denoted by 
$U_{\tau}(\supset T)$. We denote the $T$-invariant closed (fixed) 
point of $U_\tau$ as $p_\tau$. 
We write 
$\tau^\vee:=\{x\in M\otimes \R\mid \langle x,\tau\rangle 
\subset \R_{\ge 0}\},$ 
its interior as $\tau^{\vee, o}$, 
$\mathcal{S}_\tau:=\tau^\vee \cap M$, 
and $\mathcal{S}_\tau^o:=\mathcal{S}_\tau\cap 
\tau^{\vee,o}$. 

We take $(0\neq) \xi\in \tau$ which is not contained in 
$E\otimes_{\Q} \R$ for 
any proper 
$\Q$-linear subspace $E$ of $N_{\Q}:=N\otimes_{\Z} \Q$ i.e., 
of rational rank $r$. Henceforth, we simply say that $\xi$ 
is irrational. 
We introduce the monoids 
\begin{align*}
M_{\xi\ge 0}&:=\{m\in M\mid \langle m,\xi\rangle \ge 0\},\\
\Gamma_{\ge 0}&:=\langle M_{\xi \ge 0},\xi \rangle \subset \R. 
\end{align*}
The latter is not finitely generated if $r>1$. 
Here, $\langle , \rangle$ means the natural pairing which comes from the definition of $M$. Naturally, we have the following: 

\begin{lem}\label{lem:MGamma}
The natural pairing map 
$$\langle \xi, -\rangle \colon M\to \R,$$
which sends $m\in M$ to 
$\langle m,\xi\rangle =\sum_i m_i w_i$, is a group isomorphism 
onto its image $\Gamma:=\sum_i \Z w_i\subset \R$. Its restriction to 
$M_{\xi \ge 0}$ gives an isomorphism to $\Gamma_{\ge 0}$ as 
monoids. 
\end{lem}
The proof is straightforward from the fundamental theorem on 
homomorphism because the kernel vanishes due to the irrationality of 
$\xi$. 
For our later applications to relate with differential geometry, $\xi$ is often rather a given data, 
rather than $\tau$ (or even $N$). 

Also, for any algebraic scheme $V$ over $\C$, $V^{\rm an}$ (or $V(\C)^{\rm an}$) refers to the 
(complex) analytification i.e., a complex analytic space $V(\C)$ with the complex analytic (euclidean) topology. 
\end{ntt}

We prepare the bases of ``irrational deformations" in an algebraic manner, 
following \cite{Novikov1, 
Novikov2, FOOO, AJ, AI} for instance. 

\begin{defn}[Novikov-type rings and spectra]\label{def:Nov}
Under the Notation \ref{not.hr}, 
the {\it global N-ring} \footnote{N stands for Novikov as 
these rings are variants of 
so-called Novikov ring, in the contexts of Morse theory, Gromov-Witten theory, symplectic geometry, which seem to have different versions in the literatures cf., e.g., 
\cite{Novikov1, Novikov2, FOOO, AJ, AI}.} over a ring $A$ for $\xi\in N_{\R}$ is the monoid ring 
$A[\Gamma_{\ge 0}]$ and we call 
its spectrum $\Spec A[\Gamma_{\ge 0}]$ a {\it global N-spectrum}, 
which we denote as 
$\Delta_\xi^{\rm gl}(A)$ or simply 
$\Delta_\xi^{\rm gl}$ if the base ring $A$ 
is obvious from the context. 
Then, we define the following. 
\begin{enumerate}
\item 
Let $I$ be the ideal of $A[\Gamma_{\ge 0}]$  generated by the standard basis which 
corresponds to $\Gamma_{>0}:=(\Gamma_{\ge 0}\setminus 0)$,  
and define 
the {\it N-ring} over a ring $A$ for $\xi\in N_{\R}$ as 
$$A_{\xi}:=\{\sum_{i=1,2,\cdots}a_i t^{r_i} \mid a_i\in A, r_i\in \Gamma_{\ge 0}, r_1<r_2<\cdots, 
r_i\to \infty (i\to \infty) \},$$ 
again with the formal variable $t$. Note that naturally 
$A[\Gamma_{\ge 0}]\subset A_\xi.$ 
We call ${\rm Spec}(A_{\xi})$ the {\it N-spectrum} and denote it by $\Delta_{\xi}(A)$. 
Later, we often consider the case when the base $A$ is a field $\k$ and in that case we denote 
$A_{\xi}$ as $\k_{\xi}$. 
Note that it is clearly a valuation ring of rational rank $r$. 
When the base ring is obvious, we often abbreviate $\Delta_{\xi}(A)$ as $\Delta_{\xi}$. 
\item 
If $A\subset \C$ and a semigroup $\Gamma_{\ge 0}$ is a submonoid of $\R_{\ge 0}$, 
we define the 
{\it convergent N-ring} 
$A[[\Gamma_{\ge 0}]]^{\rm conv}$ 
as 
$$\{f=\sum_{i=1,2,\cdots}a_i t^{r_i} 
\in A_{\xi} 
\mid \exists \epsilon>0 \text{ s.t. }
f(t) \text{ absolutely converges for any }
t\in (0,\epsilon) 
\}.$$ 
We write its fraction field as 
$A((\Gamma_{\ge 0}))^{\rm conv},$
and denote their spectrum 
$\Spec A[[\Gamma_{\ge 0}]]^{\rm conv}$ 
(resp., 
$\Spec A((\Gamma_{\ge 0}))^{\rm conv}$) 
as $\Delta_{\Gamma_{\ge 0}}^{\rm conv}(A)$ 
(resp., 
$\Delta_{\Gamma_{\ge 0}}^{*, {\rm conv}}(A)$) 
and call them 
{\it convergent N-spectrum}. 
We often consider the case when 
$\Gamma_{\ge 0}$ is of the form 
of Notation \ref{not.hr} 
with some $\xi\in N_\R$, 
and denote them by 
$\Delta_{\xi}^{\rm conv}(A)$ 
(resp., 
$\Delta_{\xi}^{*, {\rm conv}}
(A)$). 
\end{enumerate}
\end{defn}

We can interpret the above notion 
as an affine variant of limit toric variety in 
\cite[\S 2.4]{gal}. 

\begin{prop}\label{Nsp}
We fix a base integral domain $A$ and irrational $\xi\in N_{\R}$ as Notation \ref{not.hr}. 
Consider the set of all rational polyhedral cones $\tau$ of $N_{\R}$ which include $\xi$,  
and associate the obvious partial order with respect to the inclusion. 

Then, 
the global N-spectrum 
$\Spec A[\Gamma_{\ge 0}]$ is 
canonically isomorphic, as locally ringed spaces, to 
the projective limit $\varprojlim_{\tau\ni \xi}U_{\tau}$ of the above 
projective system 
(or its sub-system which is still cofinal to it). 
Here, each $U_\tau$ is the 
affine  toric variety over $A$ 
which corresponds to $\tau$. 
If $A$ is a field $\k$, the 
N-spectra $\Delta_{\xi}(A)$ is 
the spectrum of the completion of the stalk at 
its natural prime ideal $I_\xi$ with respect to the 
valuation filtration $F^{\ge c}R_0:=\langle t^\gamma\mid \gamma\ge c  
\rangle$ (not the simple $I_\xi$-adic completion). 
Note that the valuation group is non-discrete. 
\end{prop}
\noindent
Note that each $U_{\tau}$  
over a field $\k$ has Krull dimension $r$. On the other hand, note that 
$k_\xi$ is a valuation ring whose value group is a subgroup of $\mathbb R$. 
Since every ordered subgroup of $\mathbb R$ has rank one, $k_\xi$ has Krull dimension 
$1$. 

\begin{proof}
Since the projective system 
$\{U_\tau\}_\tau$ consists of 
affine schemes, it is enough to check at the level of 
rings. Then, it is a straightforward confirmation so 
we leave the details to the readers. 
\end{proof}

In the next section, we also introduce a closely related 
algebraic stacks over $R$ (Definition \ref{stacks}.) 
More specifically, it will appear as a natural 
base of a modification of morphism 
from ${\rm Spec}(R)$ to an algebraic stack, in our 
higher analogue of semistable reduction type theorem 
(Theorem \ref{gAHLH}, Corollary \ref{gAHLHc}).

A real ($\R$) analogue of $\Delta_{\xi}$ 
is $\tilde{R}_{\xi}$ below. 

\begin{defn}[Analytic analogue]\label{realNov}
If the base ring is $A=\k=\mathbb{C}$, then we consider 
$r_{\xi}\colon \R\to N_{\R}$ which simply maps $x$ to $x\xi$. As usual, 
$\exp\colon \R\mapsto \R_{>0}\subset \C^*=\mathbb{G}_m(\C)$ 
means the exponential map $r\mapsto e^r$. 
Consider the map 
$$R_{\xi}:=({\rm id}_N\otimes \exp)\circ r_{\xi} \colon \R\to T(\mathbb{C}),$$
where ${\rm id}_N\colon N\simeq N$ is the identity map. 
We also set the notation $(\cdot)^{\xi}:=R_{\xi}\circ \log \colon \R_{>0}\to \R_{>0}\otimes N$, 
which informally means the $\xi$-th power map. 
Note that because of the irrationality of $\xi$, 
the map $(-)^{\xi}$ 
can not be extended to 
$\C^*$. We denote the universal cover of $\C^*$ as $a\colon \widetilde{\C^*}\to \C^*$, 
which is isomorphic to $\R_{>0}\times (\widetilde{U(1)}\simeq \R)$, 
and for each $\epsilon\in \R_{>0}$. 
$\widetilde{\C}_\epsilon^*$ 
denotes the preimage of 
$\{t\in \C^*\mid |t|<\epsilon\}$. 

Then, 
$R_{\xi}\circ \log\colon \R_{>0}\to T(\C)$ nevertheless extends to $\widetilde{\C^*}\to T(\C)$ as a 
holomorphic (analytic) homomorphism between complex Lie groups. 
This $\widetilde{C^*}$ is isomorphic to $\C$ as the (additive) 
complex Lie group, through the exponential map. 
We denote it as 
$$\tilde{R}_{\xi} \colon \widetilde{\C^*}\to T(\C).$$
Note that the image is not inside the real part $N\otimes \R_{>0}$ as 
it is also nontrivial along the argument direction $\widetilde{U(1)}$, 
hence $R_{\xi}\circ |\cdot|\circ a \neq \tilde{R}_{\xi}$. 
Nevertheless, they coincide at the level of absolute values i.e., 
$$(\cdot)^{\xi}\circ |\cdot|\circ a=|\cdot|^r\circ \tilde{R}_{\xi}=:|R_{\xi}|.$$
\end{defn}

A relation with the (algebraic) N-spectrum is 
partially explained by the following 
easy lemma. 

\begin{lem}\label{convan}
Under Notation \ref{not.hr} and 
Definition \ref{def:Nov} (ii), 
take any algebraic $\C$-scheme $X$. 
Then, any $\C$-morphism (as schemes) 
$\varphi\colon 
\Delta_\xi^{*,{\rm conv}}(\C)
\to X$ induces a 
holomorphic map 
$\varphi_\epsilon^{\rm an}
\colon 
\widetilde{\C}_\epsilon^*
\to X^{\rm an}$ for 
$0<\epsilon\ll 1$. 
\end{lem}

Here, $X^{\rm an}$ denotes the complex analytification i.e., 
$X(\C)$ with the complex analytic topology and the sheaf of 
analytic holomorphic functions. 

\begin{proof}
After replacing $X$ by an affine open neighbourhood of the image, write $X=\Spec(B)$ and choose finitely many generators of $B$. Their images are convergent generalized power series. Choosing a common sufficiently small radius, these series define holomorphic functions on $\widetilde{{\C}}^*_\epsilon$, where the real powers are single-valued. The algebraic relations among the generators remain valid, and hence the functions define the required holomorphic map.     
\end{proof}

Now, we generalize the notion of test configuration as follows. 
Henceforth, we stick to non-log setup i.e., 
of (often polarized) varieties 
without boundary divisors for simplicity. 

\subsection{Generalized test configurations and 
examples}

Henceforth, we always fix a base field $\k$ throughout the paper 
(which is sometimes assumed to be $\C$ if mentioned). 
The following concepts implicitly or 
explicitly (compare \cite[\S 2.2]{HanLi}, 
\cite[Appendix A3]{BJI}, \cite[\S 2.2]{Inoue}) 
appeared in the literature. We 
introduce the terms just for slight 
extensions and convenience for our later 
usage, without any originality. 

\begin{defn}\label{negdef}
\begin{enumerate}
\item 
A {\it generalized (unpolarized affine) test configuration} 
(resp., {\it generalized (polarized, projective) test configuration}) 
is the following set of data: 
\begin{enumerate}
\item free abelian groups $N$
of rank $r$. 
\item algebraic tori $T$
defined as $T=N\otimes_{\Z} \mathbb{G}_m$. 
\item a $r$-dimensional 
rational polyhedral cone $\tau (\supset \{0\})$ of $N_{\R}:=N\otimes \R$. 
 and its corresponding 
affine toric variety $U_{\tau}$ 
over $\k$. We denote the $0$-dimensional strata 
of $U_{\tau}$ as $p_\tau$ or $0=p_\tau$ for simplicity. 
\item a variety $\mathcal{Y}$ (resp., a polarized variety $(\mathcal{Y},\mathcal{L})$) 
and a $T$-action on it. 
\item faithfully flat 
$T$-equivariant affine morphism $p\colon \mathcal{Y}\twoheadrightarrow 
U_{\tau}$ 
(resp., 
$T$-equivariant projective 
flat morphism $p\colon (\mathcal{Y},\mathcal{L})\twoheadrightarrow 
U_{\tau}.$) 
\item \label{Tfppf}
In the context when general fibers are 
affine cones (in the sense of \cite{Od24a}), 
with respect to an additional algebraic $\k$-torus $T'$, 
we require more on the generalized 
(affine) test configuration as follows. 
$\mathcal{Y}$ 
admits $p$-fiberwise $T'$-action so that it can be written as 
$\Spec \mathcal{R}$ 
with a $\Gamma(U_\tau)$-algebra $\mathcal{R}$ 
which we decompose as 
$\mathcal{R}=\oplus_{m\in M}\mathcal{R}_{m}$ 
via the $T'$-action to the $T'$-eigensubmodules over 
$\Gamma(U_\tau)$, 
then each $R_{m}$ is locally free over 
$\Gamma(U_\tau)$. 
\end{enumerate}

\item 
For an affine integral reduced 
scheme $Y$ of finite type over a field $\k$ (resp., 
polarized variety $(Y,L)$ over $\k$), 
an {\it (affine) algebraic test deformation} 
(resp., {\it (polarized) algebraic test deformation} of $(T\curvearrowright Y,\xi)$), where 
$T$ is an algebraic $\k$-torus $N\otimes \G_m$ 
i.e., which has the co-character lattice $N$, 
and 
$\xi$ is an element of $N\otimes \R$, 
is a generalized (unpolarized affine) test configuration $p\colon\Y\to U_{\tau}$ 
with $p^{-1}(0)=Y$. 

Note that if $\k$ contains $\R$ and $T$ is $\k$-split, 
$\xi$ can be identified with  
its holomorphic\footnote{we take this choice rather than  real vector field, as in usual choice of 
Reeb vector field, as we want to 
work also over $\k\neq \C$}
vector field $\xi$ over the smooth locus $Y^{\rm sm}$ of $Y$. Indeed, if we consider the induced action of $T$ on $p^{-1}(0)$ and 
the natural map $a\colon N_{\R}\to H^0(p^{-1}(0),T_{p^{-1}(0)^{\rm sm}})$ induced by it 
(here, $T_{p^{-1}(0)^{\rm sm}}$ denotes the Zariski tangent sheaf of the smooth locus of $p^{-1}(0)$), 
the image ${\rm Im}(a)$ contains $a(\xi)$ 
and $a$ is often injective, so that we identify 
$\xi$ with a holomorphic vector field. 
\item 
An affine algebraic test deformation
of $(Y,\xi)$ over $U_{\tau}$ 
is 
{\it algebraic positive (resp., negative) weight deformation} 
(or $Y$ is algebraic scale up (resp., scale down) deformation) 
over $U_{\overline{\tau}}$ 
if 
\begin{enumerate}
\item 
$T\curvearrowright Y$ is a 
good action in the sense of \cite{LS13} (e.g., Fano cone as in \cite{CS, CS2}) and 
\item $-\tau \ni \xi$ (resp., $\tau\ni \xi$). 
\footnote{Our sign convention comes from \cite{Pink, CH} and 
matches with differential geometric intuition, while being opposite with 
some other algebraic literature.}
\end{enumerate}
We also propose alternative names of 
 algebraic positive (resp., negative) weight deformation: 
{\it scale-up (resp., scale-down) degeneration}, motivated by the differential geometric context. 
\end{enumerate}
\end{defn}


\begin{figure}[htbp]
\centering
\small
\begin{tabular}{@{}c@{\hspace{10mm}}c@{}}
$\vcenter{\hbox{\xy 0;<1mm,0mm>:
  (0,-12)*{};(72,70)*{};
  (8,30);(8,68) **\dir{-};
  (8,30);(65,21) **\dir{-};
  (10,64);(61,41) **\crv{(10,50)&(40,41)};
  (10,58);(57,37) **\crv{(10,45)&(37,37)};
  (10,52);(53,33) **\crv{(10,40)&(34,33)};
  (10,46);(49,30) **\crv{(10,36)&(31,30)};
  (10,31);(63,39) **\dir{--};
  (49,56)*+{\mathcal Y};
  (34,19);(34,13) **\dir{-} ?>*\dir{>};
  (38,16)*+{\scriptstyle p};
  (8,5);(59,15) **\dir{-};
  (8,5);(56,8) **\dir{-};
  (8,5);(60,-3) **\dir{-};
  (6,3)*+{p_\tau};
  (65,5)*+{U_\tau};
  (35,-10)*+{\text{scale up }(\xi\in-\tau)};
\endxy}}$
&
$\vcenter{\hbox{\xy 0;<1mm,0mm>:
  (0,-12)*{};(72,70)*{};
  (8,30);(8,68) **\dir{-};
  (8,30);(65,20) **\dir{-};
  (9,33);(18,63) **\dir{-};
  (9,32);(32,57) **\dir{-};
  (9,32);(46,50) **\dir{-};
  (9,31);(59,42) **\dir{-};
  (9,30);(62,31) **\dir{-};
  (9,30);(57,24) **\dir{-};
  (10,31);(63,38) **\dir{--};
  (50,57)*+{\mathcal Y};
  (34,19);(34,13) **\dir{-} ?>*\dir{>};
  (38,16)*+{\scriptstyle p};
  (8,5);(59,15) **\dir{-};
  (8,5);(56,8) **\dir{-};
  (8,5);(60,-3) **\dir{-};
  (6,3)*+{p_\tau};
  (65,5)*+{U_\tau};
  (35,-10)*+{\text{scale down }(\xi\in\tau)};
\endxy}}$
\end{tabular}
\caption{Scale-up and scale-down degenerations in
Definition~\ref{negdef}.}
\label{fig:scale-up-down}
\end{figure}

\begin{ex}\label{ex.tc}
Underlying family of affine test configuration (cf., \cite{CS}) is a 
generalized (unpolarized) test configuration 
with $r=1, \tau=\R_{\ge 0}$. 
Similarly, (polarized projective) test configuration (\cite{Don02}) is a 
generalized (polarized projective) test configuration 
with $r=1, \tau=\R_{\ge 0}$. 
\end{ex}

\begin{ex}\label{closure.cst}
If we consider an embedded affine variety $Y'\subset \A_{\k}^l$ and the natural 
action of $T=\G_m^l$ on the ambient space $\A_{\k}^l$ (with all weights $1$). Then we 
consider $Y'\times \{(1,\cdots,1)\} \subset \A_{\k}^l\times T$ and 
its orbit $T\cdot (Y'\times \{(1,\cdots,1)\})$. 
If the orbit closure $\overline{T\cdot (Y'\times \{(1,\cdots,1)\})}$ inside 
$\A_{\k}^l\times U_{\tau}$ is faithfully flat over $U_{\tau}$ for some (not necessarily strictly) convex 
rational 
polyhedral cone $\tau$ of $N_{\R}$, then 
$$\overline{T\cdot (Y'\times \{(1,\cdots,1)\})}\twoheadrightarrow U_{\tau}$$ 
gives a 
generalized (unpolarized) test configuration. This generalizes the classical construction of 
test configuration \`a la Mumford (\cite{GIT}) and Donaldson (\cite{Don02}) for the $r=1$ case.  
However, the required (faithful) flatness becomes nontrivial 
in general for $r>1$. 
\end{ex}

We explain more systematic constructions of 
scale-up degenerations (cf., Example \ref{loc.ex} Theorem 
\ref{loc.ex.thm}) and scale-down degenerations 
(Example \ref{ex.nwd} Lemma \ref{prep}). 

\begin{ex}[The degenerations of klt singularities to its 
K-semistable Fano cone]\label{loc.ex}
We discuss an algebraic version of the local version of the 
Donaldson-Sun theory \cite{DSII}. 
Suppose the base field $\k$ is algebraically closed of characteristic $0$. 
For any affine klt variety $X$ with its closed point $x$, there is a canonical 
quasi-monomial valuation $v_X$ with the center $x$ (\cite{Li, Blum}) and 
$W=\Spec \gr_{v_X} \mathcal{O}_{X,x}$ is a K-semistable Fano cone, 
in particular $\gr_{v_X} \mathcal{O}_{X,x}$ is finitely generated (\cite{LX, XZ}). 
Here, recall that $\gr_{v_X} \mathcal{O}_{X,x}:=\oplus_{a\in \R_{\ge 0}}\{f\in \mathcal{O}_{X,x}\mid v_X(f)\ge a\}/\{f\in \mathcal{O}_{X,x}\mid v_X(f)> a\}$. 

Choose its homogeneous generators $\overline{z}_1,\cdots,\overline{z}_l$, and we write their weights 
as $w_1,\cdots,w_l \in \R_{>0}$ in the holomorphic spectrum (cf., \cite{DSII}) 
and denote $\xi:=(w_1,\cdots,w_l)$. Naturally, we have the corresponding 
embedding $W\hookrightarrow \A_{\k}^l$. 
We then take lifts of $\overline{z}_1,\cdots,\overline{z}_l$ 
to the element of $\Gamma(\mathcal{O}_X)$, say 
$z_1,\cdots,z_l \in \Gamma(\mathcal{O}_X)$. By adding further $z_i$'s (and $\overline{z}_i$'s) 
if necessary, 
we can assume that the corresponding map 
$X\hookrightarrow \A_{\k}^l$, which maps $x$ to the origin $\vec{0}$, is an embedding. 

Recall that by a slight perturbation of $v_X$ to a divisorial valuation $v'$, 
we can easily show there is a natural isomorphism of the graded ring as 
$gr_{v_X}(\O_{X,x})=gr_{v'}(\O_{X,x})$ (cf., e.g., \cite[2.10]{LX}) 
and we take associated test configuration as the spectrum of graded Rees algebra 
(\cite[Chapter 5, section 4]{Matsumura}), 
which degenerates $x\in X$ to $W$ with the vertex. 
Note that this test configuration depends on the choice of the perturbation $v'$. 
We integrate these degenerations for various approximation $v'$ as follows. 

We apply Lemma \ref{lem:MGamma} for our weight vector $\xi=(w_1,\cdots,w_l)\in \R^l$, 
and also take a rational polyhedral cone $\tau$ of $N$ which contains $\xi$, 
to later discuss its perturbations inside $\tau$. 

We define the
 algebraic subtorus $T\subset \mathbb{G}_m^r$ as $T:=\mathbb{G}_m\otimes N$. 
 We take generators of 
the defining ideal of $X\subset \A_{\k}^l$, say $f_1,\cdots f_a$, 
which are even universal Gr\"obner basis of the defining 
ideal of $X$ (\cite[\S 5.1]{BCRV}). 

We then take small enough $\tau$ which includes $\xi$ and applies Example \ref{closure.cst} 
which results in 
a flat $T$-equivariant family $\pi_{\tau}\colon \mathcal{X}_{\tau}\twoheadrightarrow U_{\tau}$ 
whose fiber is $X$ over 
$(1,\cdots,1)$ and the fiber over the $T$-invariant point $p_{\tau}$ is $W$  
as the proof \cite[2.10]{LX} shows as $\tau$ is small enough. 

By using the fact that $f_i$s form a universal Gr\"obner basis, 
it easily follows that for any rational ray $\R_{>0}\xi'$ in $\tau$, the inducing valuation 
on $X$ is represented by $\R_{>0}\xi'$, which we denote as $v_{\xi'}$, 
its extended Rees construction (\cite[Chapter 5, section 4]{Matsumura})  
is nothing but the 
restriction of $\pi_{\tau}$ to the orbit closure $\overline{(\G_m\otimes (N\cap \R \xi'))\cdot [X]}$ 
along $\xi'$-direction. 
Further, from the Zariski openness of the irreducibility condition, by replacing $\tau$ by a smaller 
enough one if necessary, 
we can and do assume that all the fibers of $\pi_{\tau}$ are irreducible. 
(Even more is true, as \cite[2.10]{LX} again shows that the fibers over the toric boundary 
are all isomorphic.) 
Thus its inducing valuation 
on $X$ is represented by $\R_{>0}\xi'$. 
Further, from the construction, the above generalized test configuration $\pi_{\tau}$ is a 
positive weight deformation (scale up degeneration) in the sense of Definition \ref{negdef}. 

To make the construction even more canonical, to remove the ambiguity of the choice of $\tau$, 
if we take a projective limit  
with respect to the injective system of $\tau$ via inclusions, 
we obtain the canonical deformation family over $\varprojlim_{\tau} U_{\tau}$. 

We end Example \ref{loc.ex} by the following brief summary 
of our discussion. 

\begin{Thm}\label{loc.ex.thm}
For any affine klt variety $X$ with its closed point $x$, over an algebraically closed field $\k$ 
of characteristic $0$, 
consider the volume-minimizing valuation $v_X$ of $\Q$-rank $r$ in the sense of \cite{Li}, 
and $W:=\Spec ({\rm gr}_{v_X}(\mathcal{O}_{X,x}))$. 

There is an affine toric variety $U_\tau$ and an 
scale-up deformation over it $$\pi_{\tau}\colon \mathcal{X}_{\tau}\twoheadrightarrow 
U_{\tau},$$ which restricts to 
$W\times \partial U_\tau\twoheadrightarrow \partial 
U_\tau$, over the toric boundary $\partial U_{\tau}:=U_\tau \setminus T$. Its general fibers are $X$. 
\end{Thm}
See closely related studies in 
\cite[\S 2]{Teissier}  and 
\cite[the proof of 2.10]{LX} which follows {\it loc. cit.}. 

In the above, 
for the dual $N$ of the groupification $M$ of ${\rm Im}(v_X)$, 
$\tau$ is a 
small enough rational polyhedral cone in $N\otimes \R$ 
which includes $\xi=(w_1,\cdots,w_l)$ where $w_i$ are the values of $v_X$ for the lifts of homogeneous  
generators of ${\rm gr}_{v_X}(\mathcal{O}_{X,x})$. 
To remove the ambiguity of $\tau$, we can also take 
affine morphisms $\varprojlim_{\tau}\X_{\tau}\to \varprojlim_{\tau}U_\tau$, 
which is completely canonical whose general fiber is (the base change of) $X$ while the 
closed special fiber is $W$ i.e., the canonical degeneration of $X$ to $W$. In particular, if we pass to the completion of a stalk (see Proposition \ref{Nsp}), we obtain a canonical affine 
family $\X_\xi \to \Delta_\xi$ over the N-spectrum $\Delta_\xi$. 
\end{ex}

\begin{Rem}\label{ex.GIT}
As the general context of \cite{GIT}, 
consider a reductive algebraic group $G$ and its 
action on a polarized variety $(M,L)$ over an 
algebraically closed field $\k$. 
Take a semistable closed point $x$ of $X$. 

Then, the orbit closure $\overline{G\cdot x}\subset M$ 
contains a unique $G$-orbit which is closed inside 
the semistable locus. We denote it by $G\cdot y$ 
with $y\in X$. 
By \cite{Iwahori} (cf., also \cite[Chapter 2, \S 1]{GIT}), 
there is a one parameter subgroup 
$\lambda\colon \G_m \to G$ with 
$\lim_{t\to 0}\lambda(t)\cdot x=y$ 
though it is far from unique. 
If we take a maximal torus of $G$ which contains the 
image of $\lambda$, then clearly we have 
$\overline{T\cdot x}\ni y$, 
with a little less isomorphism classes. 
If $M$ is a moduli of certain 
varieties, this will give us a 
generalized test configuration but we do not 
know its canonicity in general (since 
$T$ is unique only up to the conjugation). 
\end{Rem}

\begin{Rem}\label{loc.ex2}
Let us consider isotrivial 
degenerations of K-semistable Fano cone $W$ to K-polystable Fano cone $C$, which first appeared in \cite{DSII} where the $C$ is 
the metric tangent cone of the Gromov-Hausdorff limits of 
K\"ahler-Einstein manifolds, and later \cite{LWX} gave an algebraic proof to the fact that 
such isotrivial degeneration exists for any $W$ 
as an affine test configuration. 

Recall that $C$ is given an action of algebraic split torus $T$ from the definition, and 
by the reductivity of 
$G:={\rm Aut}_T(C)$ (\cite[Appendix]{DSII}, cf., also \cite{Od24a}), which means the group of automorphisms of $C$ which commutes 
 with the $T$-action, and the action of 
 its centralizer to the multi-Hilbert scheme \cite{multiHilb} $H$ as in \cite{DSII}. 
 From these facts, 
 it is 
 easy to see similar remark apply as above 
 Remark \ref{ex.GIT}. 
\end{Rem}

\vspace{7mm}

\begin{ex}[$\R$-test configurations, K\"ahler-Ricci flow]\label{ex.r-tc}
For a general polarized variety $(X,L)$ over $\k$, 
there is a notion of $\R$-test configuration 
(cf., \cite{HanLi, BLXZ} as originated in \cite[\S 3]{CSW}) which we recall now. 
\begin{defn}
If we set $R_m:=H^0(X,L^{\otimes m})$ for $m\in \Z_{>0}$, and $R_{(X,L)}:=\oplus_{m\in \Z_{\ge 0}}
R_m$, 
the filtration $\F^{\lambda}R_m$ of $R_m$ for $\lambda\in \R$  
which satisfies the following is 
called {\it filtration} of $(X,L)$ (cf., \cite{Szek}). 
\begin{enumerate}
\item $\F^{\lambda}R_m$ is a $\k$-vector subspace of $R_m$ for each $\lambda\in \R, m\in \Z$. 
\item (decreasingness) $\F^{\lambda'}R_m\supset \F^{\lambda}R_m$ for any $\lambda'\le \lambda$. 
\item (linearly boundedness) there is $e_-$ and $e_+$ such that $\F^{e_- m}R_m=R_m$ and 
$\F^{e_+ m}R_m=0$ for any $m$. 
\item (multiplicativity) $\F^{\lambda_1}R_{m_1}\cdot \F^{\lambda_2}R_{m_2}\subset 
\F^{\lambda_1+\lambda_2}R_{m_1+m_2}$ 
for any $\lambda_1, \lambda_2 \in \R, m_1, m_2 \in \mathbb{Z}$. 
\item (as a convention: left-continuous) $\cap_{\lambda'< \lambda}\mathcal{F}^{\lambda'}R_m=
\mathcal{F}^{\lambda}R_m$. 
\end{enumerate}
For a filtration $\F^{\lambda}R_m$, 
we can define the graded ring 
\begin{align*}
{\rm gr}_{\F}(R_{(X,L)})&:=\oplus_{m\in \Z_{\ge 0}} {\rm gr}_{\F}R_m\\ 
&=\oplus_{m\in \Z_{\ge 0},\lambda\in \R} 
(\F^{\lambda}R_m/\cup_{\lambda'>\lambda} \F^{\lambda'}R_m).
\end{align*}
A filtration is called 
a {\it $\R$-test configuration} if ${\rm gr}_{\F}(R_{(X,L)})$ is finitely generated 
algebra over $\k$. 
\end{defn}

Now we briefly 
recall the following geometric realization of 
$\R$-test configuration as a generalized test configuration, 
as essentially proved in 
\cite{Teissier, HanLi, BJ, Inoue} 
under the name of 
``geometric $\R$-test configuration II" 
(\cite[p.10]{HanLi}), 
``rank $r$ test configuration" (\cite[Appendix A.3]{BJ}), 
and ``polyhedral configuration" (\cite[\S 2.2.2, Prop 2.20]{Inoue}) 
respectively. See also {\it loc. cit.} 
for more information on this correspondence. 

\begin{prop}[\cite{HanLi, BJ, Inoue}]
For a $\R$-test configuration $\F$ of $(X,L)$, with its central fiber 
$(Y:={\rm Proj}({\rm gr}_{\F}(R(X,L))),\mathcal{O}(1))$ and its natural holomorphic 
vector field $\xi$, 
we can associate 
a canonical (polarized) algebraic test deformation of 
$(Y,\mathcal{O}(1),\xi)$ over a certain affine toric variety $U_{\tau}$. 
\end{prop}

\begin{proof}
We proceed similarly as the proof of Theorem \ref{loc.ex.thm} in Example \ref{loc.ex}. 
Take homogeneous generators of ${\rm gr}_{\F}(R_{(X,L)})$ and their $\Z$-homogeneous 
lifts to $R(X,L)$ 
as $f_1,\cdots,f_l$ 
and denote their gradings as $(m_i, w_i) (i=1,\cdots,l)$ with $m_i\in \Z_{\ge 0}, w_i\in \R$. 
We set $\xi:=(w_1,\cdots,w_l)$ and take $\tau\ni \xi$ similarly. 
Then we apply Example \ref{loc.ex} to the cone $Y':={\rm Spec}R(X,L)$ with the natural 
$\G_m$-action. By dividing the resulting family by $\G_m$, we obtain the desired family. 
We omit the details as the proof is almost verbatim to Example \ref{loc.ex}. 
\end{proof}
\noindent
Similarly to the remark after Theorem \ref{loc.ex.thm}, 
if we use the language of limit toric variety $\varprojlim_{\tau} U_{\tau}$ (or its germ at the $0$-dimensional strata) in 
\cite[\S 2.4]{gal} or use its spectra of the completion (Proposition \ref{Nsp}) as the N-spectra $\k_\xi$, 
we can eliminate the 
ambiguity of $\tau$ by the set of all its subdivisions. 

Recall that \cite[Theorem 1.6]{CW} and then \cite[Corollary 1.3]{Bamler}, 
\cite[Theorem 1.1]{WZ} proved the Hamilton-Tian conjecture which says 
that the Gromov-Hausdorff limit $W$ of the normalized 
K\"ahler-Ricci flow on (possibly K-unstable) Fano manifold $X$ is $\Q$-Fano variety with 
weak K\"ahler-Ricci soliton, defined as 
$$\frac{\partial}{\partial t}
\omega_t=\omega_t-{\rm Ric}(\omega_t),$$
with $[\omega_t]\in 2\pi c_1(X).$
Then \cite{CSW} 
proved even more i.e., it is  
obtained as the central fiber of $\R$-test configuration of $(X,\O(-K_X))$, which minimizes 
the H-invariant as proved in \cite{DerSze, HanLi} (cf., also \cite{BLXZ}). 

\begin{cor}[cf., \cite{CSW, HanLi, BLXZ}]
For any 
K-unstable $\Q$-Fano variety $X$, consider the corresponding 
a (klt) 
$\Q$-Fano variety $W$ with action of an algebraic torus $T$ 
and $\xi\in ({\rm Hom}(\G_m,T)\otimes \R) \setminus 
 ({\rm Hom}(\G_m,T)\otimes \Q)$, which 
 satisfies the (modified) K-semistability in 
 the sense of \cite[\S 3.5]{BWN}, \cite[\S 2]{HanLi}. 
This is the $1$-st process of the degeneration to the limit of 
the normalized K\"ahler-Ricci flow on $X$ 
as \cite{CSW, HanLi}, 
also characterized as the linear bounded filtration 
whose $H^{\rm NA}$-invariant is minimized (\cite{HanLi, BLXZ}). 
Thus, there is a canonical algebraic test deformation of 
$W$ to $X$. 
\end{cor}

\end{ex}

\begin{ex}[Negative weight deformations]\label{ex.nwd}

Conlon-Hein \cite{CH} introduces a 
notion of negative weight deformation of 
Calabi-Yau cones 
in terms of sequence 
of affine test configurations, 
for their classification of asymptotically conical 
Calabi-Yau manifolds. We alternatively introduce 
the notion of algebraic negative deformation in 
Definition \ref{negdef} without using sequence 
(to answer Open Questions 1.4(2) of \cite{CH}). 
For our later purpose, we also clarify the relation as follows. Here, 
we assume $\k=\C$ as 
\cite{CH}. 
We take an affine $\C$-variety $Y$ and an algebraic split torus $T=N\otimes \G_m$ 
acting on $Y$, with $\xi\in N_{\R}$. 

\begin{lem}[Rational approximation - relation with \cite{CH}]\label{prep}
\begin{enumerate}
\item \label{1st}
Suppose $p\colon \mathcal{Y}
\twoheadrightarrow U_{\tau}$ 
is a scale-down degeneration (algebraic negative weight deformation) 
(Definition \ref{negdef}) 
of $(Y,\xi)$. In particular, 
$\xi\in \tau$. 
Then, its general fiber is a negative $\xi$-weight deformation in the sense of \cite[Definition 1.7]{CH}. 
\item (\cite{Od24c})\label{2nd}
Conversely, if $Y$ has only isolated singularities and $T$ acts as a good action 
in the sense of \cite{LS13}, 
any negative $\xi$-weight deformation of normal affine variety 
with $Y$ and $\xi$ 
arises in this way. 
\end{enumerate}
\end{lem}
\begin{proof}
For the item \eqref{1st}, 
we can assume that $\tau=\Z_{\ge 0}^r$, 
$N=\Z^r$. Then we approximate 
$\xi$ by 
$\xi_i=(\frac{r_{i,1}}{m_i},\frac{r_{i,2}}{m_i},\cdots,\frac{r_{i,r}}{m_i})\in 
\tau\cap N_{\Q}$, with nonnegative integers $r_{i,j} (j=1,\cdots,r)$ and 
$m_i$s so that a convergence $\xi_i \to \xi$ holds for $i\to +\infty$. 
Then take a morphism $\varphi_{\vec{r}}\colon \A^1\to U_{\tau}$ 
defined by $t\mapsto (t^{r_{i,1}},t^{r_{i,2}},\cdots,t^{r_{i,r}})$ and 
consider the base change $\mathcal{Y}\times_{p,U_\tau,\varphi_{\vec{r}}}\A^1$, 
which we denote as $p'_i\colon W'_i\to \A_1$. 
It can be naturally seen 
as an affine test configuration but we modify the $\G_m$-action $\sigma'_i$ on it 
by composing with the $m_i$-th power map $\G_m\to \G_m$. In particular, 
$p'_i(\sigma'_i(z,x))=z^{m_i}p'_i(x)$ for $z\in \G_m(\k), x\in W'_i$.  
Then, we take the base change of $p'_i\colon W'_i\to \A^1$ with $\sigma'_i$ 
by the $m_i$-th power map $q_i\colon \A^1\to \A^1$ whose target is regarded as the 
base of $p'_i$. A straightforward verification shows that the sequence, 
with respect to $i=1,2,\cdots$, of the obtained fiber product 
$p_i\colon W_i(:=W'_i\times_{p'_i,\A^1,q_i} \A^1)\to \A^1$ with the $\G_m$-action 
induced by $\sigma'_i$ satisfies the conditions of \cite[Definition 1.7]{CH}, 
with $c_i$ of {\it loc. cit.} as $m_i$, $\mu_i$ of {\it 
loc. cit.} as $1$, hence the latter being auxiliary 
parameters. 
Suppose $Z$ is 
the maximum Artin scheme 
supported at the $T$-invariant point $p_\tau \in U_\tau$, 
such that $\mathcal{Y}|_Z$ 
is trivial family. 
Such $Z$ exists due to the Zorn's lemma 
combined with the formal smoothness 
(\cite[Part I, Part IV]{EGA}) of ${\rm Aut}_T(Y)$ 
which is an algebraic group as the stabilizer of 
$[Y]$ in the multi-Hilbert scheme. 
By using this $Z$, $k_i$ 
of {\it loc. cit.} can be written as 
$m_i {\rm length}\mathcal{O}_
{\varphi_{\vec{r}}^* Z}$, 
which explains the nature of 
the $\xi$-weight in {\it loc. cit.}. 
\end{proof}
The item \eqref{2nd} is proved in \cite{Od24c} 
in which some more related materials are also explored. 
\end{ex}

\section{Higher rank semistable reduction theorem}\label{sec:3}

In this section, we generalize the theorem of \cite[\S 6]{AHLH} 
to the families over higher dimensional toroidal base 
or Novikov type rings. 

\subsection{Toroidal stack-theoretic setup}

Here are some preparations on 
stack-theoretic framework, 
following \cite{HL, Heinloth, AHLH}. 

\begin{defn}\label{stacks}
Under the notation \ref{not.hr} in the previous subsection, 
we set the following Artin stacks. 
\begin{enumerate}
\item 
Consider the quotient stack of the affine toric variety $U_{\tau}$: 
$$\Theta_{\tau}:=[U_{\tau}/T],$$
generalizing the case $r=1, \tau=\R_{\ge 0}$ 
treated in \cite{HL, AHLH}, which is denoted simply 
as $\Theta$. We denote the $T$-fixed point of $U_{\tau}$ as 
$p_{\tau}$. In \cite{HL}, 
(the isomorphic object) 
$\Theta_{\tau}$ for the case of regular simplicial $\tau$ 
    is also essentially introduced under the 
disguise of $\Theta^r$, the $r$-times self-product of $\Theta$. 
For our purpose, since we will 
encounter a setup with various $\tau$ for a fixed $N$ later, 
we prefer the above setup and 
notation. 
\item \label{stacksST}
We fix a discrete valuation ring $R$ with a valuation 
$v\colon (R\setminus 0)\to \Z_{\ge0}$. 
To set up a notion of elementary modification of 
family along $R$ (Definition \ref{def:elmod} later), 
we prepare the following. 

Further, we set $N_{\R}:=N\otimes_{\Z}\R$ and 
take a $r$-dimensional rational 
polyhedral cone $\tau'$ inside 
$N_{\R}\oplus \R$, whose generators can be taken inside $N\times 
\R_{> 0}$ and assume that the restriction of 
the first projection $p\colon N_{\R}\times \R\to N_{\R}$ 
to 
$\tau'$ 
is a PL isomorphism 
(in particular, bijectively) 
to $\tau\subset N_{\R}$. 
In other words, $\tau'$ is the graph of 
a homomorphism $\tau\to \R_{> 0}$. 

We further assume 
the normal vector of $\tau'$ can be taken as 
$(-m,1)\in ((\tau^o)\cap M)\times \Z$, 
i.e., $\langle \tau',  (-m,1)\rangle=0$. 
This is equivalent to that $p|_{\tau'}(\tau'\cap (N\oplus \Z))=\tau\cap N$. 
(Note that after replacing $R$ by its $n$-th ramified (finite) extension for 
some positive integer $n$ if necessary, 
the assumption will be satisfied.) 
Under this assumption, clearly such $m$ is uniquely and naturally determined by $\tau'$
as an element of $\tau^\vee$. 
Vice versa, $m$ determines the 
hyperplane of $M_\R \oplus \R$ spanned by $\tau'$. 
We denote the character on $T$ which corresponds to $m$ as $\chi^{m}$. 
Then we define 
\begin{align}
\label{STtilde}
\widetilde{\rm ST}_R(\tau')&:=\Spec \biggl(\dfrac{R[\tau^{\vee}\cap M][X]}
{(\pi-X\chi^{m})}\biggr),\\ 
{\overline{\rm ST}}_R(\tau')&:=[\widetilde{\rm ST}_R(\tau')/T]. 
\end{align}
\end{enumerate}
\end{defn}


\begin{figure}[htbp]
\centering
\small
\begin{tabular}{@{}c@{\hspace{4mm}}c@{}}
$\vcenter{\hbox{\xy 0;<1mm,0mm>:
  (0,0)*{};(88,76)*{};
  (72,66);(14,66) **\dir{-};
  (72,66);(18,55) **\dir{-};
  (72,66);(21,36) **[thinner]\dir{-};
  (72,66);(28,17) **\dir{-};
  (72,66);(72,5) **\dir{-};
  (69.1,66.0);(68.6,65.3) **\dir{-};
  (65.1,66.0);(64.0,64.4) **\dir{-};
  (61.1,66.0);(59.4,63.4) **\dir{-};
  (57.1,66.0);(54.7,62.5) **\dir{-};
  (53.0,66.0);(50.1,61.5) **\dir{-};
  (49.0,66.0);(45.4,60.6) **\dir{-};
  (45.0,66.0);(40.8,59.6) **\dir{-};
  (41.0,66.0);(36.1,58.7) **\dir{-};
  (37.0,66.0);(31.5,57.7) **\dir{-};
  (16.9,66.0);(15.0,63.2) **\dir{-};
  (26.7,20.5);(29.8,19.0) **\dir{-};
  (24.8,25.8);(32.5,22.0) **\dir{-};
  (22.8,31.2);(35.2,25.0) **\dir{-};
  (21.4,36.2);(37.9,28.0) **\dir{-};
  (25.4,38.6);(40.6,31.0) **\dir{-};
  (29.4,40.9);(43.3,34.0) **\dir{-};
  (33.4,43.3);(46.0,37.0) **\dir{-};
  (37.4,45.6);(48.7,40.0) **\dir{-};
  (41.4,48.0);(51.3,43.0) **\dir{-};
  (45.4,50.3);(54.0,46.0) **\dir{-};
  (49.4,52.7);(56.7,49.0) **\dir{-};
  (53.4,55.0);(59.4,52.0) **\dir{-};
  (57.4,57.4);(62.1,55.0) **\dir{-};
  (61.4,59.8);(64.8,58.0) **\dir{-};
  (65.4,62.1);(67.5,61.0) **\dir{-};
  (69.4,64.5);(70.2,64.0) **\dir{-};
  (69.6,5.7);(72.0,8.0) **\dir{-};
  (66.3,6.6);(72.0,12.3) **\dir{-};
  (63.0,7.5);(72.0,16.5) **\dir{-};
  (59.6,8.4);(72.0,20.7) **\dir{-};
  (56.3,9.3);(72.0,24.9) **\dir{-};
  (53.0,10.2);(72.0,29.2) **\dir{-};
  (49.7,11.1);(72.0,33.4) **\dir{-};
  (46.4,12.0);(72.0,37.6) **\dir{-};
  (43.1,12.9);(72.0,41.8) **\dir{-};
  (39.7,13.8);(72.0,46.1) **\dir{-};
  (36.4,14.7);(72.0,50.3) **\dir{-};
  (33.1,15.6);(72.0,54.5) **\dir{-};
  (29.8,16.5);(72.0,58.7) **\dir{-};
  (45.2,36.1);(72.0,63.0) **\dir{-};
  (30.2,19.5);(36.0,14.8) **\dir{-};
  (33.8,23.5);(49.1,11.3) **\dir{-};
  (44.6,35.5);(72.0,13.6) **\dir{-};
  (48.2,39.5);(72.0,20.5) **\dir{-};
  (51.8,43.5);(72.0,27.3) **\dir{-};
  (55.4,47.5);(72.0,34.2) **\dir{-};
  (59.0,51.5);(72.0,41.1) **\dir{-};
  (62.6,55.5);(72.0,48.0) **\dir{-};
  (66.2,59.5);(72.0,54.9) **\dir{-};
  (69.8,63.6);(72.0,61.8) **\dir{-};
  (72,66);(36,61) **\dir{-} ?>*\dir{>};
  (72,66);(44,47) **\dir{-} ?>*\dir{>};
  (44,47);(44,62) **\dir{--} ?>*\dir{>};
  (24,61)*+{\hbox{\large$\tau$}};
  (40,42)*+{\hbox{\large$\tau'$}};
  (55,22)*+{\hbox{\large$\hat{\tau}'$}};
  (49,64)*+{\ensuremath{\xi}};
  (51,50)*+{\ensuremath{\xi'}};
  (77,36)*+{\ensuremath{(\vec{0},1)}};
\endxy}}$
&
$\vcenter{\hbox{\xy 0;<1mm,0mm>:
  (0,0)*{};(90,76)*{};
  (18,58);(23,64) **[thinner]\dir{-};
  (29,58);(34,64) **[thinner]\dir{-};
  (40,58);(45,64) **[thinner]\dir{-};
  (51,58);(56,64) **[thinner]\dir{-};
  (62,58);(67,64) **[thinner]\dir{-};
  (73,58);(78,64) **[thinner]\dir{-};
  (22,47);(27,53) **[thinner]\dir{-};
  (33,47);(38,53) **[thinner]\dir{-};
  (66,47);(71,53) **[thinner]\dir{-};
  (77,47);(82,53) **[thinner]\dir{-};
  (18,36);(23,42) **[thinner]\dir{-};
  (29,36);(34,42) **[thinner]\dir{-};
  (62,36);(67,42) **[thinner]\dir{-};
  (73,36);(78,42) **[thinner]\dir{-};
  (22,25);(27,31) **[thinner]\dir{-};
  (72,26);(76,31) **[thinner]\dir{-};
  (42,14);(10,23) **\dir{-};
  (42,14);(77,25) **\dir{-};
  (42,14);(77.4,25.3) **\dir{-};
  (42,14);(79,4) **\dir{-};
  (42,14);(79.4,3.7) **\dir{-};
  (67,65);(67,14) **\dir{--} ?>*\dir{>};
  (42,14);(67,14) **\dir{--};
  (20,16)*+{\ensuremath{V(\tau')}};
  (43,11)*+{\ensuremath{p_\tau}};
  (79,15)*+{\ensuremath{U_\tau}};
  (76,45)*+{\ensuremath{{\rm Spec}(R)}};
  (50,48)*+{\ensuremath{T\times\Spec(K)}};
\endxy}}$
\end{tabular}
\caption{($r=2$) The left depicts cones $\tau, \tau', \hat{\tau}'$ and the 
right depicts 
$\widetilde{\rm ST}_R(\tau')$ (Definition \ref{stacks}).}
\label{fig:def31-schematic}
\end{figure}

The above picture is put to help human imagination, 
assuming $r=2$ just for here, and 
the (nearly) horizontal direction facing left corresponds to 
the $r=2$-dimensional cone $\tau 
(\subset N_{\R}\oplus \R_{\ge 0})$ 
and the vertical downward direction corresponds to $\R_{\ge 0} (\subset 
\tau \times \R_{\ge 0})$. Note $\hat{\tau}'$ is a $3$-dimensional cone. 

These are higher dimensional generalization of 
the case $r=1, \tau=\R_{\ge 0}$ 
introduced in \cite[\S 2.B]{Heinloth} 
and \cite[\S 3.5]{AHLH}. 
Note that, if $R$ is the stalk of $\A^1_t$ at $(t)$, 
$\biggl(\dfrac{R[\tau^{\vee}\cap M][X]}
{(\pi-X\chi^{m})}\biggr)$ 
of \eqref{STtilde} is a localization of the coordinate ring of affine toric 
variety which corresponds to 
\begin{align}
\hat{\tau}':=\tau'+\R_{\ge 0}(\vec{0},1). 
\end{align}
In general, we have a natural map 
$$b\colon {\overline{\rm ST}}_R(\tau')\to 
{\rm Spec}(R).$$ If we denote the 
fraction field of the 
valuation ring $R$ as $K$, 
the  generic fiber of $b$, which is 
a priori an open substack of ${\overline{\rm ST}}_R(\tau')$ 
is nothing but ${\rm Spec}(K)$ 
as it is straightforward to see. 
In particular, $b$ is birational in the sense that 
$b$ restricts to an isomorphism ${\rm Spec}(K)\to {\rm Spec}(K)$. 
If we denote by $\kappa$ the residue field $R/(\pi)$, 
then the special fiber of $b$ i.e., the fiber stack over the closed point 
$(\pi)\in {\rm Spec}(R)$ is a quotient stack over it: 
\begin{align}
[\biggl({\rm Spec}\dfrac{\kappa[\tau^\vee\cap M][X]}{(X\cdot \chi^m)}\biggr)/T]. 
\end{align}
It contains an irreducible component $(X=0)$ 
which is $\Theta_{\tau,\kappa}$ i.e., the $\Theta_\tau=[U_\tau/T]$ over $\kappa$. 
On the other hand, if we consider the special fiber of 
$\widetilde{\rm ST}_R(\tau')\to {\rm Spec}(R)$, then it contains 
the one dimensional strata defined by the ideal 
generated by $\pi$ and $\chi^{m}$s for all $m\in 
\tau^\vee\cap M\setminus 0$. We also denote its image 
as a closed subset of 
${\widetilde{\rm ST}}_R(\tau')$ as $V(\tau')$. 
These basic structures of $\widetilde{\rm ST}_R(\tau')$ and 
${\overline{\rm ST}}_R(\tau')$ 
will be useful in the proofs of the main theorem \ref{gAHLH} and its key 
lemma \ref{resol.indet}. 

\begin{defn}[Complex analytic analogue]
If $\k=\C$, then we can literally make sense of the toroidality of the birational morphism 
$b^{\rm an} \colon U_{\hat{\tau}'}^{\rm an}\to \C\times U_{\tau}^{\rm an}$. 
We denote the complex analytic space $(b^{\rm an})^{-1}(\Delta\times U_{\tau}^{\rm an})$ as 
$\widetilde{\rm ST}_{\Delta}^{\rm an}(\tau')$ and 
its analytic quotient stack 
${\overline{\rm ST}}_{\Delta}^{\rm an}(\tau'):=[\widetilde{\rm ST}_{\Delta}^{\rm an}(\tau')/T(\C)^{\rm an}]
$. 
\end{defn}

Now we extend the notion of $\Theta$-strata (\cite{HL}) and elementary modifications 
(\cite{AHLH}) accordingly, 
which morally corresponding to the ``optimal destabilization". The point of our generalization is that 
they are {\it not} necessarily 
one parameter algebraic families but rather higher dimensional families,  
which contains ``irrational" directions. 

Below, we fix an algebraic stack $\mathcal{M}$ over $\k$. 

\begin{defn}[Higher $\Theta$-strata]
\label{higher.theta.strata}

Consider the mapping stack ${\rm Map}(\Theta_\tau, \mathcal{M})$ 
of $\Theta_\tau$ to $\mathcal{M}$, as an algebraic stack 
by \cite[5.10, 5.11]{AHR} 
(cf., also \cite[10.17]{AHR2}, \cite[5.1]{HLP}). 
Extending \cite[Definition 2.1.1]{HL}, a 
 {\bfit{higher $\Theta$-stratum}} of type $\tau$ and 
 rank $r$ in 
$\mathcal{M}$ 
consists of a union of connected components $\mathcal{Z}^{+}\subset 
{\rm Map}(\Theta_\tau, \mathcal{M})$ such that the natural evaluation morphism 
${\rm ev}_{(1,\cdots,1)} \colon \mathcal{Z}^{+} \to \mathcal{M}$ 
at the identity $(1,\cdots,1)$ of the algebraic torus $\G_m\otimes N$ 
is a closed embedding. 
Note that, if $\tau$ is a regular simplicial one and we fix an isomorphism 
$\Theta^r\simeq \Theta_\tau$ 
for a rational polyhedral cone $\tau$, 
we can regard it as an open closed substack 
of ${\rm Filt}^r(\mathcal{M}):={\rm Map}(\Theta^r, \mathcal{M})$. 

Furthermore, if $\mathcal{M}$ is quotient stack of the form $[H/G]$ where $H$ is an 
algebraic $\k$-scheme $H$ and $G$ is an algebraic $\k$-group, 
we often write 
${\rm ev}_{(1,\cdots,1)} (\mathcal{Z}^{+})=[Z^+/G]$ for 
some 
$G$-invariant closed subscheme $Z^+\subset M$. 
\end{defn}

For the case $r=1$, the above is introduced by 
\cite[\S 2.1]{HL}. 
Now, we are ready to introduce the generalization of 
elementary modification of \cite[\S 3]{AHLH}. 

\begin{defn}[Elementary ($\tau$-)modification cf., {\cite[\S 3.5, Definition 3.36]{AHLH}} for $r=1$ case]\label{def:elmod}

An {\it extended elementary ($\tau$-)modification} of $f\colon {\rm Spec}(R)\to \mathcal{M}$ 
 is the data of a morphism 
 $\overline{f}\colon 
 \overline{\rm ST}_{R}(\tau') \to \mathcal{M}$, 
 with some $\tau'(\simeq \tau)$ as Definition \ref{stacks} 
 \eqref{stacksST}, 
 along with an isomorphism 
 $f\simeq \overline{f}|_{\rm Spec(R)}$. Here, to make sense of the restriction symbol, 
 $\rm Spec(R)$ is regarded as a localization of 
 $\overline{{\rm ST}}_R(\tau')$ since we have 
 $\Spec (R)=[\Spec R[M]/T]$. 
\end{defn}

\subsection{Higher $\Theta$-semistable reduction theorem}
\label{sec:stred}

Below, we generalize the semistable reduction type theorem in \cite[\S 6.1]{AHLH}, 
which was for $r=1$ case. 
Morally speaking, it says that a 
one parameter family which  strictly degenerates into a higher $\Theta$-strata, there is a {\it canonical}  modification of the family to 
avoid a higher $\Theta$-strata 
i.e., with ``milder degeneration". This is a version of 
semistable reduction type theorem, and later 
Theorem \ref{theta.stratification.
stablereduction} is its iterated version along a 
certain kind of stratification 
(higher $\Theta$-stratification cf., Definition \ref{Theta.strat}). 

\begin{Thm}\label{gAHLH}
Let $\k$ be an algebraically closed field of characteristic $0$,  
$(R,\mathfrak{m},v)$ 
be an equi-characteristic discrete 
valuation ring which is $\k$-algebra 
with $R/\mathfrak{m}\simeq \k$ (discrete valuation $\k$-algebra), 
and 
$\mathcal{M}$ be 
an algebraic stack over $\k$ 
which is a quotient stack $[H/G]$ where 
$H$ is a finite type scheme\footnote{henceforth we sometimes paraphrase as ``algebraic $\k$-scheme" or 
``algebraic scheme over $\k$" following the standard references, such as Mumford's Red Book 
and Stacks project.} 
over $\k$ and $G$ is an algebraic 
affine group acting on $H$. 

On the other hand, as in Notation \ref{not.hr}, we set 
$T=\mathbb{G}_m\otimes_{\Z} N$, where $N$ is a free abelian group of rank $r$, 
and fix a rational regular polyhedral cone $\tau\subset N_{\R}$ and $(0\neq )\xi\in \tau$ 
which is irrational if $r>1$. 

Consider an arbitrary morphism 
$f\colon \Spec(R)\to \mathcal{M}$, 
which lifts to $\tilde{f}\colon \Spec(R)\to H$. 
Suppose that $f(\mathfrak{m})$ sits in a {\bfit{higher 
$\Theta$-stratum}} $\mathcal{Z}^+$ in 
$\mathcal{M}$ while $f(\eta_{R})$ is not in $\mathcal{Z}^+$, where $\eta_{R}$ denotes the generic point of $R$. 
Then, after 
\begin{itemize}
\item replacing $\tau$ by a priori smaller rational polyhedral cone of $N_{\R}$ which still contains $\xi$ 
\item replacing $R$ by its finite extension 
local $\k$-morphism $R\hookrightarrow R'$, 
so that the valuation of 
$R'$ restricts to that of $R$ ($R'$ is again a DVR), 
\item 
taking a rational polyhedral cone $\tau'$ of $N_{\R}\oplus \R$ which projects as PL isomorphism 
(in particular, bijectively) 
to $\tau\subset N_{\R}$, as in Definition \ref{stacks}, 
\end{itemize}
there exists the following extension $\overline{f}$ of $f$: 

\begin{enumerate}[(i)]
\item \label{SText}
$f$ extends to its 
{\bfit{extended elementary $\tau$-modification}} (see Definition 
\ref{def:elmod}) 
\begin{align}
 F\colon {\overline{\rm ST}}_R(\tau')\to \mathcal{M}
 \end{align}
such that the generic point $\eta_{V(\tau')}$ of the 
one dimensional locus $V(\tau')$ (recall from 
Definition \ref{stacks}) 
maps to outside $\mathcal{Z}^+$. 

\item (Canonicity) \label{uniqueext}
Assume $G$ is an algebraic $\k$-torus. 
Then, $\tau'\ni \xi'=(\xi,c)$ for some $c\in \R_{>0}$ 
is unique up to taking subcone (shrinking). In particular, 
$c$ is a unique value \footnote{can be regarded as a 
certain critical value, which implicitly appears in 
\cite{Song} and was an initial inspiration, and also led 
to \cite{Od24c}. We thank S.Sun for the discussion.}
Hence, the image by $F$ of $V(\tau')$ in $|\mathcal{M}|$ 
is also unique. 
\end{enumerate}
\end{Thm}

As we mentioned, the above Theorem \ref{gAHLH} 
is applied in \cite{Od24a, Od24c} in the context of 
K\"ahler/Sasaki-Einstein geometry and corresponding 
algebro-geometric K-stability theory. 

\begin{Rem}\label{remNov}
For $F$ in the above Theorem, if we replace $\tau$ 
and $\tau'$ by smaller cones, as far as they contain 
$\xi$ and $\xi'$ respectively, one can still obtain 
$F$ as a morphism for each such smaller $\tau'$. 
By passing to the projective limit, in particular, we obtain a 
morphism from global N-spectrum to $\mathcal{M}$ (cf., Proposition \ref{Nsp}). 
This formally removes some 
redundancy of choices of $\tau'$. 
\end{Rem}

By restricting the obtained extension $\bar{f}$ of $f$, 
we obtain the semistable reduction type statement, 
generalizing \cite[\S 6]{AHLH}. 

\begin{cor}[Semistable reduction along higher $\Theta$-strata]\label{gAHLHc}
In the setup of Theorem \ref{gAHLH}, 
there exists a finite extension of DVRs $R\to R'$ and 
a modification 
$$f_{{\rm new}}\colon 
\Spec (R')\to 
\mathcal{M}$$ which maps the 
closed point outside 
$\mathcal{Z}^+$ and extends the base change of $f|_{\Spec(K)}$. 
\end{cor}

\begin{proof}[Proof of Corollary \ref{gAHLHc} (assuming 
Theorem \ref{gAHLH})]

Take a rational ray $\R_{\ge 0}\xi'$ inside 
$\tau'$ and we consider the corresponding restriction of 
$\overline{F}$ after replacing $R$ by its finite extension. 
Then, we obtain the desired map. 
\end{proof}

\begin{proof}[Proof of Theorem \ref{gAHLH}]
Our proof of Theorem \ref{gAHLH} is divided into steps. 
The core idea and arguments are in Step \ref{s1} 
with the key lemma \ref{resol.indet}, 
which may be one feature of this paper, 
and the 
rest is about how to reduce to the arguments of Step 1. 

\begin{description}
\item [{\bf Step 1}] for quotient stack of an affine scheme by an algebraic torus, 
\item [{\bf Step 2}] for general quotient stack (by 
an affine algebraic group)
\end{description}

\begin{Step}\label{s1}

First, we assume that $\mathcal{M}=[H/T]$ with an affine algebraic $\k$-scheme 
$H$ and an algebraic split $\k$-torus $T$ 
and prove the theorem in this case. Later, we essentially 
reduce to this case. Note that $f$ automatically lifts to 
$$\tilde{f}\colon \Spec(R)\to H$$ since any $T$-torsor on $\Spec(R)$ is trivial. 
From our assumption on the structure of higher 
$\Theta$-strata on $\mathcal{Z}^+$ and 
the Bia\l{}ynicki-Birula
type theorem \cite[1.4.7]{HL}, 
$\mathcal{Z}^+=[Z^+/T]$ for some closed subscheme 
$Z^+\subset H$ and 
there is a morphism $f\colon {\rm Spec}(R)\to H$ and 
$T$-equivariant morphism $Z^+\times U_{\tau}\to H$ which extends 
$f$, 
following Notation \ref{not.hr}, for some rational 
polyhedral cone $\tau$. 
Since there is also the (compatible) $T$-action on $H$, they induce a 
$T$-equivariant rational map 
\[
\Phi\colon U_{\tau}\times \Spec(R)\dashrightarrow H. 
\]
From the condition that $Z^+$ gives a union of connected components in the mapping stack, 
it follows that $\Phi$ is indeed not defined in the neighborhood of 
$(p_{\tau},\mathfrak{m})$. 
Roughly put, we partially resolve its indeterminacy in an explicit 
``(affine, non-proper) toroidal" manner 
as the following lemma. 

This is our core idea for the proof of Theorem \ref{gAHLH} and 
is of {\it purely scheme-theoretic} nature. 

\begin{lem}[key lemma]\label{resol.indet}
If necessary, by replacing $R$ by its finite totally ramified 
extension $R'$ of degree $d$ (still a DVR), 
$\tau$ by a smaller rational polyhedral cone which still contains $\xi$, 
there is some $c'\in \mathbb{R}_{>0}$ and 
$\R_{\ge 0}\xi' \subset N_{\R}\oplus \R$ 
with $\xi'=(\xi,c')$  such that the convex hull 
$$\hat{\tau}':=\tau'+\R_{\ge 0}(\vec{0},1),$$ 
is spanned by a regular rational polyhedral cone of 
$N\oplus \frac{1}{d}\Z$ and 
satisfies the following: 
\begin{enumerate}
\item there is a $T$-equivariant morphism
\[
\tilde{F}\colon \widetilde{\rm ST}_{R'}(\tau')  \longrightarrow H
\]
extending the orbit map $a\circ ({\rm id}\times f)$ on 
$T\times\Spec R'$;
\item if $\eta_{V(\tau')}(\simeq \mathbb{G}_m)$ is the generic point of the one-dimensional toric orbit
corresponding to $\tau'$, then
\[
  \tilde F\bigl(\eta_{V(\tau')}\bigr)\notin Z^+.
\]
\end{enumerate}
\end{lem}
\begin{proof}[Proof of Lemma \ref{resol.indet}]
We put $B:=R[\tau^\vee\cap M]$ and work in a $T$-invariant affine open 
$W=\Spec A\subset H$ containing $\lim_{p\to p_\tau}\Phi((p,\mathfrak{m}))$. 
This limit means the image of $p_\tau$ for the extension $U_\tau\to H$ 
of the action map 
$T\to H$ defined as $t\mapsto t\cdot \tilde{f}(\mathfrak{m})$. 
Choose homogeneous generators 
$z_1,\ldots,z_N$ of $A$.  After taking a common $T$-semi-invariant denominator, we 
write 
\[
  \Phi^*(z_i)=\frac{f_i}{g}\in{\rm Frac}(B),
  \hspace{4mm} f_i,g\in B,
\]
and set $\mathcal{J}:=(g,f_1,\ldots,f_N)$. 
We use only the standard affine open subset of its blow up 
\begin{equation}\label{eq:gchart}
  X_g:=D_+(g)\subset{\rm Bl}_{\mathcal J}(\Spec B), \hspace{3mm}
  X_g=\Spec B\!\left[\frac{f_1}{g},\ldots,\frac{f_N}{g}\right].
\end{equation}
The relations among the $z_i$ vanish in ${\rm Frac}(B)$ and therefore in the subring in 
\eqref{eq:gchart}; hence $X_g$ carries a morphism to $W\subset H$. 
Note that $D_+(g)$ is an open affine $T$-invariant subset of 
the blow up of $\mathcal{J}$ and (i) follows from the construction. 
Further, if we write $f_i, g$ as products of powers of 
$\pi$ and characters of $T$ modulo units, we obtain a 
toroidal 
description of $X_g\to U_\tau\times {\rm Spec}(R)$ 
which corresponds to subdivision of the cone 
$\tau\times \R_{\ge 0}$. 

Now we pass to a more global version of the construction 
of $X_g$ 
(without logically using $\mathcal{J}$). 
\footnote{We include both arguments because they have 
complementary
advantages: the first gives a quicker algebraic 
construction of $X_g$ 
in an explicit affine setting, while the second gives 
more geometrically visualized (hence probably more 
understandable) proof of Lemma \ref{resol.indet}.} 
For that, first we replace $H$ by the normalization of 
the 
closure of the image of ${\rm Spec}(R)\times T$.
As $H$ is affine, we can compactify to a normal projective $T$-equivariant variety 
$H\xhookrightarrow{\iota}\overline{H}$. 
Furthermore, by blowing up the boundary, 
we can and do assume that the boundary support a 
Cartier effective divisor $D$. 
We consider the normalization $\widetilde{\Gamma(\iota\circ \Phi)}$ of 
the graph $\Gamma(\iota\circ \Phi)$ 
of $\iota\circ \Phi\colon U_\tau\times {\rm Spec}(R)\dashrightarrow 
\overline{H}$ 
as $\Gamma(\iota\circ \Phi)\subset U_\tau\times {\rm Spec}(R)\times 
\overline{H}$. The projection to 
$U_{\tau}\times {\rm Spec}(R)$ composed with the normalization, 
denoted as $\widetilde{\Gamma(\iota\circ \Phi)}\to U_\tau\times {\rm Spec}(R)$ 
is a blow up of some $T$-invariant coherent ideal sheaf, 
because of the regularity of $U_\tau \times \Spec(R)$, 
which we denote as $\mathcal{J}'$. 
Then, from the construction, there is the following morphism 
from its blow up 
\begin{align}\label{maptoHbar}
{\rm Bl}_{\mathcal{J}'}(U_\tau\times {\rm Spec}(R))
\to \overline{H}.
\end{align}

The uniformizer $\pi$ of $R$ 
gives a morphism 
$c_{\pi}\colon {\rm Spec}(R)\to \A^1_{\k}$ with $c_\pi^\#(t)=\pi$, 
so that it induces an isomorphism at the formal completions. 
Since $R$ is a DVR and \(\mathcal J'\) is \(T\)-invariant, we may choose
homogeneous generators and write
\[
\mathcal J'=
(\pi^{a_1}\chi^{u_1},\ldots,\pi^{a_s}\chi^{u_s})
\subset R[\mathcal{S}_\tau].
\]
Set
\[
\mathcal J'_0=
(t^{a_1}\chi^{u_1},\ldots,t^{a_s}\chi^{u_s})
\subset \k[t][\mathcal{S}_\tau].
\]
Then
\[
\mathcal J'
=
(c_\pi\times{\rm id}_{U_\tau})^{-1}\mathcal J'_0.
\]
The ideal \(\mathcal J'_0\) is a monomial ideal for the
$(T\times\mathbb{G}_m)$-action on $U_\tau\times\mathbb A^1_k$; 
hence the normalization of its blow-up is the toric modification 
associated with a rational polyhedral subdivision $\overline{\Sigma}$ 
of 
$\tau\times\mathbb R_{\geq0}.$ We denote its fiber-product by 
$$\mathcal U_{\overline{\Sigma}}:=T_{N\oplus \Z}{\rm emb}(\overline{\Sigma})\times_{\A^1_{\k}}\Spec(R).$$
Now we restrict the preimage of $H$ in it via the morphism 
\eqref{maptoHbar}. 
Note that the complement $\overline{H}\setminus H$ 
supports the pullback of Cartier divisor $D$, 
which is $T$-invariant. We denote it by $E$. 
Further, the dense subset 
\[
  T\times \Spec(K)
\]
maps into $H$, $E$ is supported on
the toroidal boundary of $\mathcal U_{\Sigma}$. 
Thus, 
there is a sub-fan $\Sigma'$ of $\overline{\Sigma}$ which is 
still $r+1$-dimensional i.e., with nontrivial interior. 
We emphasize that it does not necessarily 
contain $(\xi,0)$. 
In this way, we obtain 
\[
  \tilde{F}\colon\mathcal U_{\Sigma'}\longrightarrow H. 
\]

(As a convenient phrasing, in the continuation of our proof, 
we call these $\mathcal{U}$-type schemes e.g., 
$\mathcal U_{\Sigma'}$ or $\mathcal U_{\overline{\Sigma}}$ toroidal  (models).) 
Further, if necessary, after shrinking $\tau$ to a smaller rational polyhedral regular 
cone which still contains $\xi$ and taking the intersection of 
$\overline{\Sigma}$ with $\tau\times \R_{\ge 0}$, we 
may arrange so 
that the intersection with the interval between  $(\xi,0)$ and 
$(0,1)$ is divided by a chain of $r$-dimensional walls. Let 
\[
  \tau'_1,\ldots,\tau'_\ell
\]
be the walls, ordered from the vertical side to the horizontal side,
and put $C_i:=\overline{V(\tau'_i)}\simeq\mathbb{P}^1 (\subset 
\mathcal{U}_{\overline{\Sigma}})$. 
Among these walls, we in particular are interested in 
those whose corresponding orbits are 
contained in the toric open subset 
$\mathcal U_{\Sigma'}=\overline F^{-1}(H)$. 
The toroidal strata corresponding to chamber on the 
vertical side maps to $\tilde{f}(\mathfrak m)\in Z^+$, whereas on 
the toroidal strata  chamber on the 
horizontal side, either $\tilde{F}$ is not well-defined or 
maps outside $Z^+$. Indeed, otherwise the retraction to
$V(\tau\times\{0\})\simeq\Spec R$, together with the fact that $\mathcal Z^+$ is a
union of connected components of ${\rm Map}(\Theta_\tau,\mathcal M)$, would force
$f(\eta_R)\in\mathcal Z^+$, where $\eta_R$ 
is the generic point of $\Spec(R)$. 

Let $j$ be the smallest index such that
\begin{equation}\label{eq:first-outside}
  \widetilde F\bigl(\eta_{C_j}\bigr)\notin Z^+, 
\end{equation}
where $\eta_{C_j}$ is the generic point of $C_j$. 
Such $j$ exists by the previous paragraph, 
and for every positive integer 
$i<j$ the closedness of $Z^+$ gives 
$\tilde{F}(\eta_{C_i})\in Z^+$.  We claim that $\tilde{F}|_{C_i}$ is 
constant and prove it below. 


\begin{figure}[htbp]
\centering
\small
\begin{tabular}{@{}c@{\hspace{1mm}}c@{}}

\hspace{-6mm}
\raisebox{6mm}{
$\vcenter{\hbox{\xy 0;<1mm,0mm>:
  (0,0)*{};(86,78)*{};
  (72,62);(18,62) **\dir{-};
  (72,62);(20,40) **\dir{-};
  (72,62);(28,24) **\dir{-};
  (72,62);(40,12) **\dir{-};
  (72,62);(57,10) **\dir{-};
  (72,62);(72,8) **\dir{-};
  (44,68)*+{\ensuremath{\tau\subset N_{\mathbb R}}};
  (16,36)*+{\ensuremath{\tau_j}};
  (25,18)*+{\ensuremath{\tau_{j-1}}};
  (36,23)*+{\ensuremath{\cdots}};
  (40,8)*+{\ensuremath{\tau_i}};
  (52.5,21)*+{\ensuremath{\cdots}};
  (58,6)*+{\ensuremath{\tau_1}};
  (78,18)*+{\ensuremath{(\vec{0},1)}};
\endxy}}$}
&
\hspace{-3mm}
\raisebox{-2mm}{$\vcenter{\hbox{\xy 0;<1mm,0mm>:
  (0,0)*{};(106,82)*{};
  (12,63);(22,45) **\dir{-};
  (22,45)*{\bullet};
  (22,45);(37,32) **\dir{-};
  (22.8,45.4);(37.8,32.4) **\dir{-};
  (37,32);(41,29) **\dir{.};
  (41,29);(55,23) **\dir{-};
  (41.8,29.4);(55.8,23.4) **\dir{-};
  (55,23);(60,21.7) **\dir{.};
  (60,21.7);(72,20) **\dir{-};
  (60.8,22.1);(72.8,20.4) **\dir{-};
  (72,20)*{\bullet};
  (72,20);(90,20) **\dir{-};
  (72,20.4);(90,20.4) **\dir{-};
  (14,50)*+{\ensuremath{C_j}};
  (27,36)*+{\ensuremath{C_{j-1}}};
  (45,36)*+{\ensuremath{\cdots}};
  (46,21)*+{\ensuremath{C_i}};
  (65,17)*+{\ensuremath{C_1}};
    (60,30)*+{\ensuremath{\cdots}};
  (74,17)*+{\ensuremath{p_\tau}};
  (85,17)*+{\ensuremath{U_\tau}};
  (80,79);(23,55) **\dir{--};
  (78,75);(34,41) **\dir{--};
  (79,73);(51,28) **\dir{--};
  (81,73);(68,26) **\dir{--};
  (83,78);(83,26) **\dir{--};
  (92,43)*+{\ensuremath{{\rm Spec}(R)}};
  (23,33);(23,30) **\dir{-};
  (23,30);(47,10) **\dir{-};
  (47,10);(81,10) **\dir{-};
  (81,10);(81,12) **\dir{-};
  (54,7)*+{\ensuremath{\in Z^+}};
\endxy}}$}
\end{tabular}
\caption{Walls and curves (rank $r=1$ case) in the proof of Lemma 3.8.}
\label{fig:walls-in-the-proof}
\end{figure}

 
As a preparation, we write the $T$-equivariant morphism as follows (cf., \cite[1.4.1, 1.4.8]{HL}): 
\begin{equation}\label{eq:universal}
  \mu\colon U_\tau\times Z^+\longrightarrow H,
  \qquad \mu(1,z)=z. 
\end{equation}

Let $g_i\colon C_i\to Z^+$ be the induced morphism.  The wall $\tau'_i$ is the
graph of a nonzero nonnegative rational linear form on $\tau$; hence, after
multiplication by a positive integer, we may choose $\lambda_i\in\tau^\circ\cap N$
whose corresponding $\G_m$-action on the open orbit of $C_i$ is nontrivial and contracts it, as $s\to0$,
to a toric fixed point $q_i\in C_i$. 
If we denote the associated $\G_m$-action to $\lambda_i$ 
simply as $\lambda_i$, then it 
extends to a morphism 
$\A^1\to U_\tau$, and \eqref{eq:universal} gives
\[
  H_i\colon\A^1\times C_i\longrightarrow H,
  \hspace{1cm} H_i(s,x)=\mu(\lambda_i(s),g_i(x)).
\]
For $s\neq0$, equivariance gives
\[
  H_i(s,x)=\lambda_i(s)\cdot\tilde{F}(x)
          =\tilde{F}(\lambda_i(s)\cdot x).
\]
It follows that $H_i(0,x)=\tilde{F}(q_i)$ on the dense torus orbit of $C_i$, and
hence on all of $C_i$ because $H$ is separated.  The rigidity lemma applied to
$H_i\colon\A^1\times C_i(\simeq \mathbb{P}^1)\to H$ now shows that 
$H_i|_{\{1\}\times C_i}=\tilde{F}|_{C_i}$ is constant 
as desired. 

Let $\rho_i\colon\mathcal U_{\Sigma'}\to\mathcal U_{\Sigma'(i)}$ be 
the toroidal 
contraction obtained by removing the wall $\tau'_i$.  Since $\rho_i$ is proper, its
target is normal, and $(\rho_i)_*\mathcal O_{\mathcal U_{\Sigma'}}=
\mathcal O_{\mathcal U_{\Sigma'(i)}}$, the morphism $\tilde{F}$ descends through 
$\rho_i$ whenever it is constant on $C_i$.  Indeed, near the contracted point, choose
an affine neighbourhood of $\tilde{F}(C_i)$ in $H$; properness allows us to
shrink the base so that its inverse image maps into this affine open, and the above
equality of structure sheaves gives the descended morphism.

Contracting successively $C_1,\ldots,C_{j-1}$ 
(in other words, replaces $\Sigma'$ by the corresponding coarsening) 
makes $\tau'_j$ adjacent to the 
vertical ray. Then, 
\[
  \widehat\tau'=\tau'_j+\R_{\geq0}(0,1)
\]
is a cone of the resulting fan.  We put $\tau':=\tau'_j$. 


\begin{figure}[htbp]
\centering
\small
\begin{tabular}{@{}c@{\hspace{1mm}}c@{\hspace{1mm}}c@{}}
\raisebox{-1mm}{$\vcenter{\hbox{\xy 0;<0.68mm,0mm>:
  (0,0)*{};(106,82)*{};
  (12,63);(22,45) **\dir{-};
  (22,45)*{\bullet};
  (22,45);(37,32) **\dir{-};
  (22.8,45.4);(37.8,32.4) **\dir{-};
  (37,32);(41,29) **\dir{.};
  (41,29);(55,23) **\dir{-};
  (41.8,29.4);(55.8,23.4) **\dir{-};
  (55,23);(60,21.7) **\dir{.};
  (60,21.7);(72,20) **\dir{-};
  (60.8,22.1);(72.8,20.4) **\dir{-};
  (72,20)*{\bullet};
  (72,20);(90,20) **\dir{-};
  (72,20.4);(90,20.4) **\dir{-};
  (14,50)*+{\ensuremath{C_j}};
  (27,36)*+{\ensuremath{C_{j-1}}};
  (45,36)*+{\ensuremath{\cdots}};
  (46,21)*+{\ensuremath{C_i}};
  (60,30)*+{\ensuremath{\cdots}};
  (65,17)*+{\ensuremath{C_1}};
  (74,17)*+{\ensuremath{p_\tau}};
  (85,17)*+{\ensuremath{U_\tau}};
  (80,79);(23,55) **\dir{--};
  (78,75);(34,41) **\dir{--};
  (79,73);(51,28) **\dir{--};
  (81,73);(68,26) **\dir{--};
  (83,78);(83,26) **\dir{--};
  (94,43)*+{\ensuremath{{\rm Spec}(R)}};
  (23,33);(23,30) **\dir{-};
  (23,30);(47,10) **\dir{-};
  (47,10);(81,10) **\dir{-};
  (81,10);(81,12) **\dir{-};
  (54,7)*+{\ensuremath{\in Z^+}};
\endxy}}$}
&
$\vcenter{\hbox{\xy 0;<0.75mm,0mm>:
  (0,-15)*{};(36,18)*{};
  (4,0);(32,0) **\dir{-} ?>*\dir{>};
  (18,-9)*+{\ensuremath{\substack{\text{contracting}\\ C_1,\ldots,C_{j-1}}}};
\endxy}}$
&
\raisebox{-1mm}{$\vcenter{\hbox{\xy 0;<0.68mm,0mm>:
  (0,0)*{};(88,82)*{};
  (13,57);(35,25) **\dir{-};
  (35,25)*{\bullet};
  (35,25);(79,25) **\dir{-};
  (35,25.4);(79,25.4) **\dir{-};
  (72,72);(25,46) **\dir{--};
  (72,72);(72,29) **\dir{--};
  (45,67)*+{\ensuremath{\tau'=\tau_j}};
  (12,42)*+{\ensuremath{V(\tau')}};
  (38,18)*+{\ensuremath{p_\tau}};
  (71,18)*+{\ensuremath{U_\tau}};
  (83,48)*+{\ensuremath{{\rm Spec}(R)}};
\endxy}}$}
\end{tabular}
\caption{contracting $C_1,\ldots,C_{j-1}$, 
which corresponds to removal of walls $\tau'_1,\cdots,
\tau'_{j-1}$.}
\label{fig:contract-in-the-proof}
\end{figure}


We write the simplicial subcone $\tau'$ as 
$\sum_{i=1}^r \mathbb{R}_{\ge 0}v'_i$, 
where $\{v'_i\}_i \cup (\vec{0},1)$ are generators of the extremal rays of $\hat{\tau}'$. 
Write $v'_i=(v_i,\frac{l_i}{d_i})$ with $l_i,d_i\in \mathbb{Z}_{>0}$ and set 
$d:={\rm lcm}\{d_i\}$. Then, $v'_i\in N\oplus \frac{1}{d}\Z$ and 
$v'_1,\ldots,v'_r,\left(0,\frac{1}{d}\right)$ 
form its basis and $\hat{\tau}'$ is (spanned by) 
its regular rational polyhedral cone. 

Then, we take an extension of $R$ as follows; we consider 
$R\to R'$ with ramification degree $d$, which is obtained as a localization 
of $R[x]/(x^d-\pi)$ at $(x)$. 

Note $\widetilde{\rm ST}_{R'}(\tau')$ is well-defined 
due to the regularity of $\hat{\tau}'$, 
and is identified with the normalization of 
$U_{\hat{\tau}'}\times_{\A^1_{\k}} R'$. 
Further, 
(base change of) $\tilde{F}$ is extendable to 
$\widetilde{\rm ST}_{R'}(\tau')\to H$ 
from the construction, 
and 
${\overline{\rm ST}}_{R'}(\tau')\to \mathcal{M}$ is obtained as its descent. The desired properties are satisfied from 
the construction and hence we finish the proof of Lemma \ref{resol.indet} \eqref{SText}. 
Note that the extension $R\to R'$ does not change when we 
shrink $\tau'$, as it follows from the construction. 

Now we can complete the proof of Lemma \ref{resol.indet}. 
Indeed, restricting the descended morphism to
$U_{\widehat\tau'}$ proves Lemma \ref{resol.indet} (i), while 
\eqref{eq:first-outside} proves Lemma \ref{resol.indet} (ii). 
\end{proof}

Note that the above Lemma \ref{resol.indet} 
implies (i) of Theorem \ref{gAHLH} for the situation of Step 1 
i.e., when $\mathcal{M}=[H/T]$.

Now we show the remaining assertion for Step 1 case 
i.e., canonicity (uniqueness) assertion \eqref{uniqueext} of Theorem \ref{gAHLH}. 
Let 
\[
  (R'_a,\tau'_a,\widetilde F_a),\hspace{1cm} a=1,2,
\]
be two choices obtained from Lemma \ref{resol.indet}, and write
\[
  \xi'_a=(\xi,c_a)\in{\rm relint}(\tau'_a).
\]
After replacing $R'_1$ and $R'_2$ by a common finite totally
ramified extension, shrinking $\tau$ around $\xi$, and passing to
a common regular refinement of the two toroidal models, we may compare
the two morphisms on one toroidal model. Denote the resulting common
morphism on the common $H$-valued toroidal open subset by
\[\widetilde{F}^{\dagger}\colon 
\mathcal{U}_{\Sigma^\dagger}\to H
\]
with $\Sigma^{\dagger}\ni \tau'_1, \tau'_2$. 
The two pullbacks agree there, since
they extend the same orbit map on the dense torus 
$T\times \Spec(K)$ and $H$ is separated. 

We recall the characterization of the wall selected in the proof of 
Lemma \ref{resol.indet}.  If $q_a\in V(\tau'_a)$ is its 
general closed point, 
then
\begin{align}\label{*****}
  \widetilde F_a(q_a)\notin Z^+,
\end{align}
whereas the orbit map
\[
  T\longrightarrow H,\hspace{1cm} t\longmapsto t\cdot q_a,
\]
extends to $U_{-\tau}\to H$, and the fixed point
$p_{-\tau}\in U_{-\tau}$ is mapped to $Z^+$. 
We call this image the $U_{-\tau}$-limit point. 
As a consequence of these for $a=1$ case, 
with the convention of 
the proof of Lemma \ref{resol.indet}, we observe the following: 
after inserting any further rational wall by subdivision, a wall 
$\tau^{(c)}\in \Sigma^\dagger$ through $(\xi,c)$ with $c>c_1$ satisfies
\[
  \tilde{F}^{\dagger}\bigl(\eta_{V(\tau^{(c)})}\bigr)\in Z^+,
\]
whereas, if $c<c_1$, it occurs on the other side of the first-exit
wall and the $U_{-\tau}$-limit of a closed point of 
$V(\tau^{(c)})$ lies outside $Z^+$.  This description is unchanged by
subdivision of fan (in $N_\Q\oplus \Q$): any additional wall 
which contains $(\xi,c)$ with $c>c_1$ 
has image in 
$Z^+$ and the corresponding rational curve in 
$\mathcal{U}_{\Sigma^\dagger}$ 
is contracted by the rigidity argument in the proof of 
Lemma \ref{resol.indet}. 

Thus, if $c_2>c_1$, the first (\eqref{*****} above) 
of the two defining properties fails 
for $\tau'_2$ i.e., 
$\tilde{F}_2(q_2)\in Z^+$. 
By the same argument applies with $a=2$, if 
$c_2<c_1$.  Hence
$c_1=c_2$. 
Note that for each positive real number $c$, if $(\xi,c)$ is contained in 
a rational hypersurface, then it is defined only by $c$ because $\xi$ is assumed to be 
irrational if $r>1$. Thus, the neighborhood of $\xi'$ in $\tau'_1$ 
and $\tau'_2$ coincides i.e., lies in the same rational hyperplane 
of $N\oplus \R$. This proves \eqref{uniqueext} of Theorem 
\ref{gAHLH} and in particular Step 1 completes. \end{Step}


\begin{Step}

In the second step, we consider the general quotient stack case i.e., 
 when $\mathcal{M}$ is a quotient stack of the form 
$[H/G]$ where $H$  is an algebraic scheme over $\k$ acted by an affine 
algebraic $\k$-group $G$. 
We essentially 
reduce the proof of the case to quotient stacks by algebraic tori. 

Let us denote $\mathcal{M}:=[H/G]$, where 
$H$ is a finite type $\k$-scheme 
acted by an algebraic $\k$-group $G$. Take an 
algebraic torus $T$ as previous notations. 
Additionally, we write the multiplication as  
$m\colon T\times T\to T$ 
and the action as $a\colon G\times M\to M$. 
We introduce the following concepts for the reduction to 
quotient stack by an algebraic split torus.

\begin{Property}[Local liftability]\label{lift}

Let
\[
\mathcal{Z}^+
\subset
\underline{\rm Map}(\Theta_\tau,\mathcal{M})
\]
be a higher $\Theta_\tau$-stratum, and write
\[
{\rm ev}_1\colon
\mathcal{Z}^+
\xrightarrow{\ \sim\ }[Z^+/G]
\hookrightarrow
\mathcal{M}
\]
for its evaluation, where
$1\in T\subset U_\tau$ denotes the identity. 
Through the above, we identify $\mathcal{Z}^+$ as a closed substack 
of $\mathcal{M}$. 

\begin{enumerate}

\item\label{lift1}
For a point $s\in |\mathcal{Z}^+|$, we say that
$\mathcal{Z}^+$ is
{\it smoothly locally liftable at $s$}
if there exist

\begin{itemize}

\item an affine algebraic $\k$-scheme $H'$ with an algebraic
$T$-action, and a $T$-invariant closed subscheme
\[
H''\hookrightarrow H',
\]

\item a smooth morphism
\[
p\colon [H'/T]\longrightarrow\mathcal{M}
\]
whose image contains $s$,

\end{itemize}

such that the following conditions hold:

\begin{enumerate}

\item\label{action}
there is an isomorphism
\[
[H''/T]
\xrightarrow{\ \sim\ }
[H'/T]\times_{\mathcal{M}}\mathcal{Z}^+
\]
which is compatible with the natural closed immersion
$[H''/T]\hookrightarrow[H'/T]$;

\item\label{limit}
the $T$-action on $H''$ extends to an algebraic monoid action
\[
U_\tau\times H''\longrightarrow H'',
\]
and the morphism
\[
[H''/T]
\longrightarrow
\underline{\rm Map}(\Theta_\tau,\mathcal{M})
\]
induced by this $U_\tau$-action and by
$p\colon[H'/T]\to\mathcal{M}$
factors through
\[
\mathcal{Z}^+
\subset
\underline{\rm Map}(\Theta_\tau,\mathcal{M}).
\]

\end{enumerate}

\item\label{lift2}
We say that $\mathcal{Z}^+$ is
{\it smoothly locally liftable}
if it is smoothly locally liftable at every point
$s\in|\mathcal{Z}^+|$ in the sense of \eqref{lift1}.

\end{enumerate}

Note that we do not assume that
\[
[H''/T]\subset[H'/T]
\]
itself gives a higher $\Theta_\tau$-stratum of $[H'/T]$.
Instead, we require the Cartesian identity in \eqref{action}
and the compatibility with the original higher
$\Theta_\tau$-stratum $\mathcal{Z}^+$ in \eqref{limit}.

\end{Property}

In practice, many 
examples of (higher) $\Theta$-strata satisfies even 
above Zariski local liftability condition \eqref{lift1}. 
The main point of the above property is that it enables us the proof of Theorem \ref{gAHLH} to the same arguments as Step \ref{s1} through the covering quotient stack by algebraic tori. 

Now we prove the following. 

\begin{Claim}\label{quot.lift}
Take a quotient stack $[X/G]$ over $\k$ and 
any higher $\Theta$-stratum $\mathcal{Z}^+=[Z^+/G]$ of type $\tau$, 
rank $r$. 
Then, it is at least smoothly locally liftable 
in the sense of 
Property \ref{lift} \eqref{lift2}. 
\end{Claim}

\begin{proof}[proof of Claim \ref{quot.lift}]

By applying \cite[4.4.6]{HL} (cf., also \cite[3.2]{AHLH}), 
there is an affine scheme $V$, an algebraic torus $T'$ 
and a smooth representable surjective morphism 
$q\colon [V/T']\to \mathcal{M}$ such that 
${\rm Filt}^r([V/T'])\to {\rm Filt}^r(\mathcal{M})$ 
is again smooth representable surjection for each $r\in 
\Z_{>0}$. 

From the Bia\l{}ynicki-Birula decomposition 
(\cite[4.3]{BB}, \cite[Theorem 1.4.7, 1.4.8]{HL}, 
cf., e.g., also \cite{Dri, DriGai}), 
the preimage of $\mathcal{Z}^+(\subset {\rm Filt}^r(\mathcal{M}))$ 
to ${\rm Filt}^r([V/T'])$ 
can be written as $[Z'/T']$ where $Z'$ is a finite disjoint union of 
affine connected components of $V^{\psi,+}$ for a homomorphism 
$\psi\colon T\to T'$. We obtain 
$$f_0\colon [Z'/T]\to \mathcal{Z}^+,$$ which is smooth 
since ${\rm char}(\k)=0$. We now apply \cite[2.8]{AHLH} 
to $\mathcal{Z}^+\subset \mathcal{M}$ with $f_0$, 
so that we obtain an affine algebraic $\k$-scheme $H'$ 
with $T$-action together with a smooth $\k$-morphism 
\begin{align}\label{goodcover}
    p\colon [H'/T]\to \mathcal{M} 
\end{align}
such that $[H'/T]\times_{\mathcal{M}}\mathcal{Z}^+\simeq [Z'/T]$. 
We complete the proof of the Claim. 
\end{proof}

We put $[H''/T]=p^{-1}(\mathcal{Z}^+)$ with $H''\subset H'$. 
Using the above Claim \ref{quot.lift}, 
now we can follow the arguments of Step \ref{s1}. 
We take $[H'/T]\to \mathcal{M}$ as Claim \ref{quot.lift} and 
 a lift $g\colon {\rm Spec}(R)\to [H'/T]$ after 
an \'etale extension of the DVR $R$, using the smoothness of $p$. 
Further, we lift it to 
\[\tilde{g}\colon \Spec(R)\to H',\]
which is possible since any $T$-torsor on $\Spec(R)$ is trivial. 
By the construction, $\tilde{g}(\mathfrak{m})\in H''$, 
$\tilde{g}(\eta_R)\notin H''$. We now apply the same 
arguments of Lemma \ref{resol.indet} to $g$ with 
composing with $p$ and use the original strata $\mathcal{Z}^+$. 
In this way, we obtain an extended elementary modification 
whose distinguished point lies outside $H''$ (which is the 
preimage of $\mathcal{Z}^+$ to $H'$). Composing with 
$p$, we obtain the required modification with closed point sent outside 
$\mathcal{Z}^+$. We complete the proof of Step 2.

\end{Step}
\end{proof}


We write the following corollary to the above 
proof of Theorem \ref{gAHLH}, 
which is about how the semistable reduction of 
\cite{AHLH} changes along the variation of direction $\xi$. 

For that, we prepare the following, to connect to 
the original $\Theta$-strata (\cite{AHLH}). 

\begin{defn}[Restriction to (usual $r=1$) $\Theta$-strata]
Consider $\xi\in \tau\subset N_{\R}$, 
partly following Notation \ref{not.hr}, 
but now we take $\xi$ to be integral 
i.e., $\xi\in N\cap \tau$ (until \eqref{AHLH.deform}). 
If we set $\iota\colon \Z \xi\hookrightarrow N$ 
the natural inclusion, note that 
there is a 
$({\rm exp}\otimes \iota) 
\colon 
\G_m \to T$ 
(compare Definition \ref{realNov}). 

For a higher $\Theta$-stratum 
of the cone $\tau$ on an algebraic 
stack $\mathcal{M}$, which we denote as 
$\mathcal{Z}^{+}\subset 
{\rm Map}(\Theta_\tau, \mathcal{M})$, 
consider restriction map 
$({\rm exp}\otimes \iota)^* \colon 
\mathcal{Z} 
\to 
{\rm Map}(\Theta, \mathcal{M})$ and 
suppose it gives the $\Theta$-strata structure 
for each $\xi$ i.e., the images are connected. 
Then we call the 
obtained $\Theta$-strata as 
{\it restriction along $\xi$} of the 
$\Theta$-strata $\mathcal{Z}$ and 
often denote as $\mathcal{Z}_\xi$. 
\end{defn}
 
\begin{cor}[Effect of variation of the direction]
\label{AHLH.deform}
In the setup of Step~1 of Theorem \ref{gAHLH}, fix an irrational
$\xi_0\in {\rm relint}(\tau)$. Suppose that the restriction of
$\mathcal{Z}^+$ along each rational ray considered below is an
ordinary $\Theta$-stratum. After shrinking $\tau$ around $\xi_0$,
there is a finite rational polyhedral subdivision $\Sigma_\tau$ of
$\tau$ with the following property. For
$\lambda\in N_{\Q}\cap{\rm relint}(C)$, let
$f_{{\rm new},\lambda}$ be the semistable reduction obtained by
restricting the modification of Theorem \ref{gAHLH} along
$\R_{\geq 0}\lambda$. Then
\[
 [f_{{\rm new},\lambda_1}(0)]
 =
 [f_{{\rm new},\lambda_2}(0)]
 \quad\text{in }|\mathcal{M}|
\]
whenever $\lambda_1,\lambda_2\in N_{\Q}\cap{\rm relint}(C)$ for
some $C\in\Sigma_\tau$. 
In particular, this point is locally constant 
for rational directions sufficiently close to $\xi_0$. 
\end{cor}

\begin{Rem}
(The above Corollary \ref{AHLH.deform} 
is a revised exposition of 
Cor 3.13 of the original arXiv v1.) 
A closely related result is also discussed in 
\cite{BHLINK} for application to generalize Donaldson-
Thomas 
type theory in a derived setup \cite{BDINKP}. 
Compare \cite[Constancy theorem, \S 6.1 and \S 7.3 case 
(i)]{BHLINK}.
\end{Rem}

\begin{proof}
Use the common toroidal model
$\widetilde{F}\colon U_{\Sigma'}\to H$ constructed in the proof of
Lemma \ref{resol.indet}. Since $\Sigma'$ is finite, after shrinking
and subdividing $\tau$, both the ordered sequence of walls met by
the segment joining $(\lambda,0)$ to $(\vec{0},1)$ and the first wall
whose generic orbit is mapped outside $Z^+$ are constant for
$\lambda\in{\rm relint}(C)$. Hence all such rational-ray restrictions
meet the same $T$-orbit in the corresponding locus $V(\tau'_C)$.
The $T$-equivariance of $\widetilde{F}$ shows that their special
points have the same image in $|[H/T]|=|\mathcal{M}|$.
Theorem \ref{gAHLH}(ii) shows that this image is independent of the
chosen common refinement. Finally, $\xi_0$ lies on no rational wall,
since it has rational rank $r$, which proves the last assertion.
\end{proof}

It often occurs in various moduli problems 
that a priori we only have 
a uniform degeneration along $\xi$-direction 
but not along a higher dimensional toric variety 
$U_\tau$. Nevertheless, we have the following remark 
after Proposition \ref{Nsp}. In the sequel papers 
\cite{Od24a, Od24c}, we use similar ideas. 

\begin{Rem}\label{limlim}
The natural inverse system of toric varieties $\{U_\tau\}_\tau$ 
provides finite-type approximations to the non-finite-type source
\(\Delta_\xi^{\mathrm{gl}}\) (cf., Proposition \ref{Nsp}). 
For a finite-type affine $T$-scheme $Y$, every
$T$-equivariant morphism
$$
  B\times\Delta_\xi^{\mathrm{gl}}\longrightarrow Y
$$
with an affine $\k$-scheme $B$ factors through
$B\times U_\tau$ for some $\tau\ni\xi$. 
Indeed, after writing $Y=\Spec(A)$ and choosing finitely
many $T$-homogeneous generators of $A$, the image of
these generators involves only finitely many exponents
of $\Gamma_{\geq 0}$. These exponents are contained in
$\tau^{\vee}\cap M$ for some rational polyhedral cone
$\tau$ containing $\xi$. This gives the required
factorization through $B\times U_{\tau}$. 
 
A closely related but more systematic 
construction of real and real-cone filtration 
stacks was subsequently given in \cite[\S7.2]{BHLINK} after our paper. 
They define ${\rm Filt}_K(\mathcal X)$ and discuss as a certain 
mapping stack. 
We do not use this construction in the proof of our main theorem. 
\end{Rem}

Our analytic (but somewhat weak) 
version of the semistable reduction theorem \ref{gAHLH} 
is as follows. In particular, we only discuss quotient stacks 
by {\it algebraic tori}, mainly due to the lack of enough references on 
analytic stacks (but cf., e.g., 
\cite[Appendix A]{Inoues} and 
references therein). 

\begin{cor}[Analytic 
analogue of Theorem \ref{gAHLH} and Corollary \ref{gAHLHc}]\label{gAHLH.R}
Let 
$\Delta:=\{t\in \C\mid |t|<1\}$ and 
$H$ an affine algebraic $\C$-scheme equipped with an algebraic 
action by an algebraic torus $T$ of rank $r$. 

Consider an arbitrary holomorphic map $\tilde{f}\colon \Delta \to H^{\rm 
an}$ induces 
$f\colon \Delta \to \mathcal{M}:=[H^{\rm  an}/T^{\rm an}]$ 
whose target is regarded as a 
complex analytic stack. 
Suppose that for a (algebraic)  
higher $\Theta$-stratum $\mathcal{Z}^+$ of type $\tau$ 
with ${\rm ev}_{(1,\cdots,1)}(\mathcal{Z}^+)=Z^+\subset H$, 
$\tilde{f}^{-1}((Z^+)^{\rm an})\subset \Delta$ 
is set theoretically $\{0\}\subset \Delta$. 
We fix an irrational element $\xi$ i.e., 
with the rational rank $r$ 
in $\tau\subset N_{\R}$ as in Notation \ref{not.hr}. 

Then, after the following preparation; 
\begin{itemize}
\item For some $m\in \mathbb{Z}_{>0}$, 
we replace $\Delta$ by $\Delta':=\{s\in \C\mid |s|<\epsilon\}$ with 
$\epsilon\in (0,1]$ and 
$\Delta'\to \Delta$ which maps 
$s$ to $t=s^d$ for some $d\in \Z_{>0}$, 
\item replace $\tau$ by a smaller rational polyhedral cone in $N_{\mathbb R}$ containing $\xi$, 
\item 
take a rational polyhedral cone $\tau'$ of $N\oplus \Z$ which projects as PL isomorphism 
(in particular, bijectively) 
to $\tau\subset N_{\R}$, as in Definition \ref{stacks}, 
\end{itemize}
$\tilde{f}$ extends to a $T(\C)^{\rm an}$-equivariant morphism 
${\widetilde{\rm ST}}^{\rm an}_{\Delta'}(\tau')\to H$, 
hence 
extends $f$ to 
$\overline{f}\colon {\overline{\rm ST}}^{\rm an}_{\Delta'}(\tau')\to \mathcal{M}$, 
and maps $V(\tau')$ 
outside $Z^+$. 
\end{cor}

\begin{proof}
Write $H=\Spec(A)$.
The germ of $\widetilde f$ at the origin induces a
$\mathbb{C}$-algebra homomorphism
$A\to \mathbb{C}[[t]]^{\rm conv}.$ 
Thus it defines an algebraic morphism 
\[
\tilde{f}_R\colon
\Spec(R)\longrightarrow H, \hspace{3mm}
R:=\mathbb{C}[[t]]^{\rm conv}. 
\]
See Definition 
\ref{def:Nov} (ii), which we apply with $\Gamma_{\ge 0}=\Z_{\ge 0}$. 
The assumption on
$\widetilde f^{-1}((Z^{+})^{\rm an})$
implies that the closed point of $\Spec(R)$ is mapped
to $Z^{+}$, whereas its generic point is not.

Thus, we may apply Theorem \ref{gAHLH} Step 1 
with 
$R=\C[[t]]^{\rm conv}$, so that 
Corollary \ref{gAHLHc} implies the desired assertion. 
Indeed, note that the finite totally ramified extension occurring in the
proof of Lemma \ref{resol.indet} may be taken in the form
\[
R=\mathbb{C}[[t]]^{\rm conv}
\longrightarrow
R'=\mathbb{C}[[s]]^{\rm conv},
\hspace{1cm}
t\mapsto s^d,
\]
with some $d\in \Z_{>0}$ if necessary. 
As a result, we obtain a 
$T^{\rm an}$-equivariant morphism 
\[
\widetilde{\rm ST}_{\Delta'}^{\rm an}(\tau')
\longrightarrow H^{\rm an}, 
\]
with some $\tau'\ni (\xi,c)$ for some $c\in \R_{>0}$. 
This completes the proof. 
\end{proof}


\subsection{Semistable reduction 
along higher $\Theta$-stratification: 
by repetition of Corollary  \ref{gAHLHc}}
We use the following setup, which is a straightforward extension of the notion of 
$\Theta$-stratification of \cite[\S 2.1]{HL}. 

\begin{defn}[Higher $\Theta$-stratification]
\label{Theta.strat}
Let $\mathcal{M}$ be an algebraic $\k$-stack with affine
diagonal.
A finite higher $\Theta$-stratification of $\mathcal{M}$
consists of the following data.

\begin{enumerate}
\item
A finite totally ordered set $\Gamma$ with maximal
element $1$.

\item
Closed substacks
\[
\mathcal{M}_{\leq c}\subset\mathcal{M},
\hspace{1cm} c\in\Gamma,
\]
such that
\[
\mathcal{M}_{\leq c'}
\subset
\mathcal{M}_{\leq c}
\quad\text{for }c'<c,
\hspace{1cm}
\mathcal{M}_{\leq 1}=\mathcal{M}.
\]

\item
For each $c\in\Gamma$, set
\[
\mathcal{M}_{<c}
:=
\bigcup_{c'<c}\mathcal{M}_{\leq c'},
\hspace{1cm}
\mathcal{M}_{\geq c}
:=
\mathcal{M}\setminus\mathcal{M}_{<c},
\]
and
\[
\mathcal{M}_{=c}
:=
\mathcal{M}_{\leq c}\cap\mathcal{M}_{\geq c}.
\]
We require that $\mathcal{M}_{=c}$ be the evaluation
image of a higher $\Theta$-stratum of some type
$(N_c,\tau_c)$ in the open substack
$\mathcal{M}_{\geq c}$. That is, we consider 
the mapping stacks ${\rm Map}(\Theta_{\tau_{c}}, 
\mathcal{M}_{\ge c})$ and 
\[
{\rm ev}_1\colon
\mathcal{Z}_c
\hookrightarrow
\mathcal{M}_{\geq c},
\hspace{1cm}
{\rm ev}_1(\mathcal{Z}_c)=\mathcal{M}_{=c}.
\]
\end{enumerate}
\end{defn}

By repeating Theorem \ref{gAHLH} and Corollary \ref{gAHLHc} 
for each higher $\Theta$-strata, we obtain their following 
version. 

\begin{Thm}[Semistable reduction along 
higher $\Theta$-{\it stratification}]
\label{theta.stratification.
stablereduction}
Following the notation of Definition \ref{Theta.strat}, we take an 
algebraic stack $\mathcal{M}$ which is a quotient stack $[H/G]$ over $\k$, 
and its higher $\Theta$-stratification 
$\{\mathcal{M}_{<c}\}_c$ 
as Definition \ref{Theta.strat}. 

Consider any $\k$-morphism $f\colon \Spec(R)\to \mathcal{M}$ for a 
discrete valuation $\k$-ring 
$(R,\mathfrak{m})$ with its fraction field $K$, which satisfies that 
$f(\Spec(K))\in \mathcal{M}_{=1}$ while $f(\mathfrak{m})\in \mathcal{M}_{< 1}$. 
Then, there is a finite sequence of modifications of $f$ and $R$  which we denote as 
$$f_i\colon  \Spec(R_i)\to \mathcal{M} (i=0,1,\cdots,l)$$ 
with $R\to R_i$ 
such that 
\begin{itemize}
\item each $(R_i,\mathfrak{m}_i)$ is a discrete valuation 
ring with injective 
local morphisms
$\k\to R=R_0 \subset R_1\subset \cdots \subset R_l$, 
\item $f=f_0$, 
\item $f|_{\Spec(K)}\times_{\Spec(K)}\Spec(K_i)=f_i|
_{\Spec(K_i)}$ 
for each $i$, 
\item $f_i(\mathfrak{m}_i)\in \mathcal{M}_{=c_i}$ with $c_0<c_1<\cdots<c_l=1$, 
\item in particular, $f_l(\mathfrak{m}_l)\in \mathcal{M}_{=1}$. 
\end{itemize}
For each $f_i$, if we 
consider the source stack $\Theta_i\simeq \Theta^{r_i}$ for their ranks  
$r_i$ and 
fix an 
isomorphism to $\Theta_{\tau_i}$ for some $N_i, \tau_i$, 
together with choices $\xi_i\in \tau_i\setminus N_\Q$. 
If $G$ is an algebraic torus, the above sequence of modifications $f_i$ can be 
taken canonically. 
\end{Thm}

\begin{proof}

We apply (the proof of) 
Theorem \ref{gAHLH} to $f$ with the strata 
$\mathcal{M}_{={c_0}}$ and take the restriction of the extension obtained $\overline{f}$ to the closed substack corresponding to a ray 
$\R_{\ge 0}v$ in $v\in N_{\Q}\cap \tau'$ in the proof of Theorem \ref{gAHLH}, 
and denote it by $f_1$. 
Then, we iterate the process. In each step, the index strictly increases so that 
it terminates after finite times because the index set is finite. 
The last statement about the canonical choices 
of $f_i$ follows also from Theorem \ref{gAHLH}. We 
can characterize such $f_i$ by the natural condition 
in terms of $\xi_i$ by Theorem \ref{gAHLH} but we avoid the 
lengthy repetition. 
\end{proof} 

Note that from the proof and canonicity statement of Theorem \ref{gAHLH}, 
$f_i$ are taken almost canonically in the sense that $f_i(\mathfrak{m}_i)$ are unique, 
and have certain characterizations. 

The author applies the above results effectively used in the 
problems around K-stability, K\"ahler geometry and moduli spaces in separate papers (cf., 
\cite{Od24a, Od24c}). 



\subsubsection*{Acknowledgement}
This paper was written in 2024. 
The author thanks E.Inoue and M.Jonsson 
for the comments on the manuscript. 
During the research, the author was partially 
supported by Grant-in-Aid for Scientific Research (B) 23H01069, 21H00973, 
Grant-in-Aid for Scientific Research (A) 21H04429, 
20H00112, 
and Fund for the Promotion of Joint International Research (Fostering Joint International Research) 23KK0249. 
Part of this work is done during 
the author's stay at Institut de math\'ematiques de Jussieu - Paris Rive Gauche. 
We appreciate the kind hospitality. 
While revising the paper from earlier versions, 
we restricted the situation and 
cut many technically involved parts for simplification e.g., 
by assuming the valuation ring is DVR over characteristic $0$. 
Furthermore, we employed generative AI (non-agent) to 
refine Claim 1 and its proof in Step 2 of Theorem 
\ref{gAHLH}. It was also used to significantly improve the 
exposition - 
such as adding Figures 1 through 3, correcting 
minor errors, and polishing English phrasing—all of which 
were thoroughly checked by the authors. We emphasize that 
the mathematical core, including definitions, key ideas, 
and methods (e.g., Step 1 of the proof of Theorem 
\ref{gAHLH}), remains largely our own. 

\subsection*{Notes added (after arXiv v3): about \cite{BHLINK}}

About $8$ and a half months after the first 
appearance on arXiv and submission of this paper, 
we learnt 
about preparation of 
\cite[\S 7.3]{BHLINK} which somewhat overlaps our \S 3.2. 
We would like to express our gratitude to 
D.~Halpern-Leistner and T.~Kinjo for informing us of his work prior to its appearance on arXiv and, above all, for the invaluable and fruitful discussions on the possible  connections and fascinating differences between our approaches. 
These discussions have been insightful, and we sincerely appreciate the time and effort. While there are notable overlaps of our papers in certain aspects, our methods and intended applications are quite distinct at the same time. We hope that our works complement each other and contribute meaningfully to future research.

Below, we provide a list of comments on the differences and overlaps between our works, reflecting our discussions with Tasuki Kinjo and also Chenjing Bu, 
for which I am truly grateful. First, we discuss the overlaps and then 
we illustrate differences. 

\begin{enumerate}

\item 
The basic objects in the two papers are closely related; 
${\rm Filt}^\Sigma(\mathcal{X})$ in their paper and 
${\rm Map}(\Theta_\tau,\mathcal{M})$ in our paper 
(with different notation: $\Sigma=\tau$, $\mathcal{X}=\mathcal{M}$). 
Generally, \cite{BHLINK} 
takes more general stacky viewpoint while we take a ``minimalistic" 
or at least more 
elementary (scheme theoretic) approach for applications. 
For instance, we directly work on ${\rm Map}(\Theta_\tau,\mathcal{M})$ 
and the corresponding notion of higher $\Theta$-stratifications while {\it loc. cit.}
works on a 
``real $\Theta$-strata (or -stratification)" 
of \cite[\S 7.3.2, \S 7.3.3]{BHLINK} by their categorical reformulation, 
which in particular allows a finite base change of $T$. 
(Their ``real filtration" is closely related to our Remark \ref{limlim}, 
in which we later deleted some part of our original discussions.)

\item \label{diff2}
Their 
 \cite[Theorem 7.3.5]{BHLINK} substantially overlaps 
our 
Corollary \ref{gAHLHc} (also Theorem \ref{gAHLH}, 
Theorem \ref{theta.stratification.
stablereduction}) despite a technical difference, and gives a 
different proof. 
Most notably, their proof of 
\cite[Theorem 7.3.5]{BHLINK} 
crucially relies on the earlier original theorem 
\cite[Theorem 6.5]{AHLH}, while our proof of Corollary \ref{gAHLHc} (also Theorem \ref{gAHLH}, 
Theorem \ref{theta.stratification.
stablereduction}) does not, so that it 
also gives alternative proof to 
\cite[Theorem 6.5]{AHLH} in the case of 
quotient stacks. 

\item 
Also, there are various differences even at the level of statements -  our statements of 
Corollary \ref{gAHLHc} (and Theorem 
\ref{theta.stratification.
stablereduction}) 
and \cite[Theorem 7.3.5]{BHLINK}. 
In our 
Corollary \ref{gAHLHc}, unlike \cite[7.3.5]{BHLINK}, 
we do {\it not} need the 
assumption on $\tau$ to be small enough (see 
the assumption on {\it loc. cit.} (ii)). 
On our side, a certain smallness condition on the cone is only required for the existence statement of the extension $\overline{f}$ part (i.e., the stronger statements) 
in Theorem \ref{gAHLH}.  

In our semistable reduction statements in  Corollary  \ref{gAHLHc}, a smaller cone does not appear and is not required. 

\item (In the original versions of our Theorem \ref{gAHLH} and 
Corollary \ref{gAHLHc}, we discuss families over 
more general valuation rings, while \cite[\S 6]{AHLH} and \cite[\S 7.3]{BHLINK} work over DVRs. However, we eliminated this 
extension for simpler exposition in the current version.) 

\item 
This seems to be the most crucial difference, 
at least from the perspective of the author. 
Our proofs and approaches to the aforementioned similar results - namely, our Corollary \ref{gAHLHc} (and Theorem 
\ref{theta.stratification.
stablereduction}) and \cite[Theorem 7.3.5]{BHLINK}
seem to be {\it rather different}. 

Indeed, we use more scheme theoretic techniques, 
notably Lemma \ref{resol.indet}, 
unlike \cite{BHLINK}, we do not rely on the original 
semistable reduction theorem \cite[Theorem 6.3, 6.5]{AHLH}, 
rather reprove them. 

On the other hand, 
the proof of  \cite[7.3.5]{BHLINK} uses 
more stack theoretic arguments and reduces to 
the original \cite[Theorem 6.3, 6.5]{AHLH}. 

\item 
More precisely, \cite{BHLINK} has two main theorems 
both under certain assumptions (cf., \eqref{qcgp..}) 
on the algebraic stacks - 
a certain local constancy of components of 
the filtration stack and the graded points stack 
in their \S 6.1 (Theorem 6.1.2), 
and finiteness of 
``special faces" in their \S 6.2 (Theorem 6.2.3, cf., also 6.2.5 and 6.2.8). 
Then, in the last subsection \S 7.3, 
they reduce the semistable reduction theorem similar to ours 
to \cite[6.3, 6.5]{AHLH} 
by applying those two theorems, in particular. 
Our Corollary \ref{AHLH.deform} is also 
parallel to their two main theorems 
in a certain aspect.

\item \label{qcgp..}
We assume $\mathcal{M}$ is global quotient stack in our theorems, 
while both \cite[6.1.2 resp., 6.2.3]{BHLINK} assumes a somewhat  technical conditions on  
algebraic stacks that is to have 
``finite cotangent weights" 
(resp., to have only ``quasi-compact graded points"). 
These are comparable to our work as quotient stacks 
are good examples which satisfy both conditions 
See {\it op. cit.} \S 5.3.3 and 6.2.2 
for the explanation and also 
\cite[Example \S 5.3.5]{BHLINK} for non-quotient algebraic 
stack  which do not satisfy their assumptions. 

\item 
The discussions in \cite{BHLINK}, notably its 
\S 5.2, rely on 
the theory of 
derived stacks and 
the cotangent complex, 
which we do not use. 

\end{enumerate}


\vspace{5mm} \footnotesize \noindent
Contact: {\tt odaka@kurims.kyoto-u.ac.jp} \\
Research Institute for Mathematical Sciences, Kyoto University, Kyoto 606-8502. JAPAN \\

\end{document}